\documentclass[12pt]{amsart}  
\usepackage{amsmath, amscd}  

\DeclareMathOperator{\rad}{\mathrm{rad}} 
\DeclareMathOperator{\op}{\mathrm{op}} 
 
\DeclareMathOperator{\init}{\mathrm{in}}

\DeclareMathOperator{\cone}{\mathrm{cone}}

\DeclareMathOperator{\conv}{\mathrm{conv}}

\DeclareMathOperator{\sort}{sort} 
\DeclareMathOperator{\sd}{sd} 
\DeclareMathOperator{\esd}{esd}

\DeclareMathOperator{\Hilb}{Hilb}

\newcommand{\toric}{toric }

\newcommand{\colim}[1]{\mathop{\underset{\underset{#1}{\longrightarrow}}{ \text{\rm lim}}}}
\newcommand{\prolim}[1]{\mathop{\underset{\underset{#1}{\longleftarrow}}{ \text{\rm lim}}}}

\newcommand{\Z}{\mathbb{Z}}  
\newcommand{\R}{\mathbb{R}}  
\newcommand{\N}{\mathbb{N}}  
\newcommand{\Q}{\mathbb{Q}}  
  
\newcommand{\Gc}{\mathcal{G}}  
  
\newcommand{\Ic}{\mathcal{I}}  
\newcommand{\Inew}{I}
\newcommand{\Pc}{\mathcal{P}}  
\DeclareMathOperator{\id}{id}  
  
\DeclareMathOperator{\ini}{in}

\DeclareMathOperator{\Ker}{Ker}

\DeclareMathOperator{\projdim}{proj\,dim}

\DeclareMathOperator{\reg}{reg}

\DeclareMathOperator{\Spec}{Spec}  
  
\DeclareMathOperator{\supp}{supp}

\DeclareMathOperator{\Tor}{Tor}   
\DeclareMathOperator{\Gen}{Gen}   

\DeclareMathOperator{\pnt}{\raise 0.5mm \hbox{\large\bf.}}

\newtheorem{thm}{\bf Theorem}[section]   
\newtheorem{lem}[thm]{\bf Lemma}  
\newtheorem{cor}[thm]{\bf Corollary}  
\newtheorem{prop}[thm]{\bf Proposition}

\theoremstyle{definition}
\newtheorem{defn}[thm]{\bf Definition}  
  
\newtheorem{rem}[thm]{\bf Remark}  
\newtheorem{ex}[thm]{\bf Example}  

\title{Subdivisions of Toric complexes}  
  
\author{Morten Brun} 
\address{FB Mathematik/Informatik, Universit\"at Osnabr\"uck, 49069 Osnabr\"uck, Germany}
\email{brun@mathematik.uni-osnabrueck.de} 

\author{Tim R\"omer} 
\address{FB Mathematik/Informatik, Universit\"at Osnabr\"uck, 49069 Osnabr\"uck, Germany}
\email{troemer@mathematik.uni-osnabrueck.de}

\begin{document}  
  
\begin{abstract}  
We introduce toric complexes as polyhedral complexes 
consisting of rational cones together with a set of integral generators for
each cone, and we define their associated face rings. 
Abstract simplicial complexes and rational fans can be
considered as toric complexes, and the face ring for toric complexes
extends Stanley and Reisner's face ring for abstract simplicial
complexes \cite{Stanley_book} and Stanley's face ring for rational
fans \cite{ST87}.   
Given a toric complex 
with defining ideal $I$ for the face ring we give a geometrical interpretation 
of the initial ideals of $I$ with respect to weight orders
in terms of subdivisions of the toric complex
generalizing a theorem of Sturmfels in \cite{ST}.
We apply our results 
to study edgewise subdivisions of abstract simplicial complexes.
\end{abstract}

\maketitle 

\section{Introduction}
The aim of this paper is to define and study toric complexes. An
embedded toric complex is a rational fan together with a distinguished
set of generators, consisting of lattice points, for each cone. Since an
abstract simplicial complex with $d$ vertices corresponds to a
simplicial fan whose rays 
are spanned by the elements of the standard basis for
$\R^d$, an abstract simplicial complex can be considered as an example
of a toric complex.  
Generalizing the Stanley-Reisner ideal of a
simplicial complex and the toric ideal of a configuration of lattice
points, we establish a connection
between binomial ideals in a polynomial ring and toric complexes. Our
main result shows that a regular
subdivision of a toric complex corresponds to the radical of an
initial ideal of the associated binomial ideal. As an application of
this correspondence we consider edgewise subdivisions of simplicial
complexes and we describe the associated deformations of Veronese
subrings of 
Stanley-Reisner rings.

Let us go more into detail.
$K$ denotes a noetherian commutative ring. 
A finite subset $F$ of $\Z^d$ determines  
a (non-normal) affine toric variety $X_F = \Spec(K[M_F])$.
Here $M_F$ is the submonoid of $\Z^d$ generated by $F$ and 
$K[M_F]$ is the subring of $K[\Z^d] =
K[t_1,\dots,t_d,t_1^{-1},\dots,t_d^{-1}]$ generated by
monomials of the form 
$t^a = \prod_{i=1}^d t_i^{a_i}$ for $a = (a_1,\dots,a_d) \in F$.
In the classical situation, where $M_F$ is the set of lattice points in
the set $\cone(F) \subseteq \R^d$ of positive real linear combinations
of elements of $F$,
the variety $X_F$ is a normal affine toric 
variety.
A {\em regular subdivision} of $\cone(F)$ supported on $F$, 
called projective subdivision in \cite[p. 111]{Mumford}, 
is a rational fan of the form $\Sigma = \{ \cone(G)
\colon G \in \Pi\}$ for a set $\Pi$ of subsets of $F$ with $\cone(F)
= \cup_{G \in \Pi} \cone(G)$ satisfying that there exist
linear forms $\alpha_G 
\colon \R^d \to \R$ such that the restrictions $\alpha_G
\colon \cone(G) \to \R$ assemble to a continuous convex function $f
\colon \cone(F) \to \R$. As explained for instance in
\cite[pp. 28--30]{Mumford} such regular subdivisions correspond to
coherent sheaves on $X_F$, and there exists a regular
subdivision of $\cone(F)$ such that the fan $\Sigma$ defines a
resolution of singularities  
$X(\Sigma) \to X_F$.

Our main result in Section \ref{subdivisions}
generalizes the correspondence established by Sturmfels in
\cite{ST_Tohoku} and \cite{ST}
between regular subdivisions of 
$\cone(F)$  
and initial ideals of the kernel $\Ic_F$ of the   
homomorphism  from the
polynomial ring $K[F] = K[x_a \colon \, a \in F]$ to $K[\Z^d]$ taking
$x_a$ to $t^a$. 

More precisely, 
a function $\omega \colon F \to \R$ defines a
weight order on $K[F]$
with 
$\omega(x^u) = \sum_{a \in F} u(a)\omega(a)$ 
for a monomial $x^u = \prod_{a \in F} x_a^{u(a)}$.  
The {\em initial ideal} $\ini_\omega(\Ic_F)$ 
of $\Ic_F$ with respect to
$\omega$ is the ideal generated by the {\em initial polynomials}
$\ini_\omega(f) = \sum_{\omega(u) = \omega(f)}
f_u x^u$ 
where $f  = \sum_{u} f_u x^u$ is non-zero in
$K[F]$ and 
$\omega(f) = \max \{\omega(u) \in \R \colon f_u \ne 0 \}$.

The function $\omega$ also defines a
collection $\Pi_{\sd_{\omega}T(F)}$ of subsets of $F$. Given $G \subseteq F$ let $C$
denote the minimal face of $\cone(F)$ containing $G$ and let $G'
= F \cap C$. The subset $G$ of $F$ is in $\Pi_{\sd_{\omega}T(F)}$ if and only if there
exists a linear form 
$\alpha_G$ on $\R^d$ with  
$\alpha_G(a) \le \omega(a)$ for $a \in G'$ such that equality holds
precisely if $a \in G$. 
This collection $\Pi = \Pi_{\sd_{\omega}T(F)}$ of subsets of $\Z^d$ is an 
{\em embedded toric complex} $\sd_\omega T(F)$, that is, 
it satisfies firstly that for $H\subseteq G$ with $G \in \Pi$  
we have that $H$ is in $\Pi$ 
if and only if 
$\cone(H)$ is a face of $\cone(G)$ 
with $\cone(H) \cap G = H$ and secondly that if $G$ and $H$ are
in $\Pi$, then their intersection is an element of $\Pi$. 
The fan $\Sigma = \{ \cone(G) \colon G \in \Pi \}$ is a {\em regular partial
subdivision} of $\cone(F)$ in the sense that 
the restrictions $\alpha_G \colon \cone(G) \to \R$ assemble to a
continuous function $f \colon \cup_{G \in \Pi} \cone(G) \to \R$ whose
restriction to any face of $\cone(F)$  contained in $\cup_{G \in \Pi}
\cone(G)$ is convex.

We carry the above idea one step further. Given an embedded toric
complex $T$ consisting of the collection $\Pi_T$ of subsets of $\Z^d$ and
a
function of the form $\omega 
\colon \Gen(T) = \cup_{G \in \Pi_T} G \to 
\R$, we have the embedded toric complex $\sd_{\omega} T(F)$ for every $F 
\in \Pi_T$. The embedded toric complex $\sd_{\omega}T$ is defined to consist
of 
the collection $\Pi_{\sd_{\omega}T} =
\bigcup_{F \in \Pi_T} 
\Pi_{\sd_{\omega} T(F)}$. This is the {\em regular partial subdivision of
  $T$} induced by $\omega$.

In analogy with the Stanley-Reisner ideal of an abstract simplicial
complex 
(see \cite{Stanley_book} for details) there is a
square-free monomial ideal
$J_T$ in $K[\Gen(T)]$ such that a monomial $x^u = \prod_{a \in \Gen(T)}
x_a^{u(a)}$ is {\em not} in $J_T$ if and 
only 
if there 
exists $F \in \Pi_T$ such that $\supp (x^u) \subseteq F$.
Here $\supp(x^u)$ is the support of $x^u$,
that is, 
the set of elements $a \in \Gen(T)$ with $u(a) \ne 0$.
Let $I_F \subset K[\Gen(T)]$ denote the ideal generated by the image of  
$\Ic_F \subset K[F]$ 
under the inclusion $K[F] \subseteq K[\Gen(T)]$ 
for $F \in \Pi_T$ and let $I_T = \sum_{F \in \Pi_T} I_F + J_T$.
As observed in 
\cite[Proposition 4.8]{EIST} 
the {\em face ring}
$K[T] = K[\Gen(T)]/I_T$ of the 
embedded toric complex $T$ agrees with the face ring of the rational
fan $\{\cone(F) \colon F \in \Pi_T\}$ considered in
\cite[Section 4]{ST87} in the case where $M_F = \cone(F) \cap 
\Z^d$ for every $F \in \Pi_T$. 
Restricted to embedded toric complexes, the statement of our main
result in 
Section \ref{subdivisions} is:
  \begin{thm}
If $\Gen(\sd_{\omega T}) = \Gen(T)$ and $\cup_{F \in \Pi_T} \cone(F) =
\cup_{G \in \Pi_{\sd_{\omega} T}} \cone(G)$, then
the radical ideal $\rad(\ini_\omega(I_T))$ of
$\ini_\omega(I_T)$ agrees with the ideal $I_{\sd_{\omega T}}$.
\end{thm}
The {\em face ring} $K[T]$ of an arbitrary  toric complex $T$, as defined in
\ref{face_ring},
is a $K$-algebra of the form $K[X]/I_T$
where $K[X]$ is a polynomial ring and 
$I_T$ is a  binomial ideal, that is, an ideal generated by
binomials of the form 
$x^u - x^v$ and by monomials $x^u$.
In the case where $K$ is a field 
every  binomial
ideal $I$ in $K[X]$
gives rise to a toric complex $T(I)$. In general it will not be possible to
represent the toric complex $T(I)$ as an embedded toric complex. Theorem
\ref{char_face} contains a precise criterion on $I$ classifying the ideals for
which the face ring
$K[T(I)]$ is 
isomorphic to $K[X]/I$.

Simplicial toric complexes as defined in \ref{define_properties} correspond
to abstract simplicial 
complexes in a strong sense. In particular, there is an abstract simplicial
complex 
associated to a simplicial toric complex, 
and there is a simplicial toric complex associated to every abstract simplicial
complex. 
The face ring of a simplicial
toric complex agrees with the Stanley-Reisner ring of the associated
abstract simplicial complex. Thus the ideals associated to simplicial toric
complexes are monomial, so on the level of ideals, simplicial toric complexes
correspond to the binomial ideals that happen to be monomial. It is hence very
restrictive to require  
a toric complex to be simplicial. 

In order to subdivide a toric complex $T$ we must have enough generators to be
able to obtain additional faces. If $T$ is simplicial, there are not enough
generators to do this. The solution to this problem is to
add generators in a controlled way.
In \ref{define_multiple} we present a
process of 
multiplying a toric 
complex by a natural number, and for $r \ge 2$ the multiple $rT$ of a toric
complex $T$ has enough generators to 
possesses subdivisions. 
On the level of face rings, multiples of toric complexes
roughly correspond to Veronese subrings of graded $K$-algebras, and
their subdivisions correspond to deformations of the Veronese ring.  
The underlying fan of a multiple $rT$ of an embedded toric complex $T$ agrees
with the underlying fan of $T$ itself, and the underlying fan of a subdivision
of $rT$ is a subdivision of the underlying fan of $T$.

For the simplicial toric complex $T(\Delta)$  constructed from an abstract
simplicial complex $\Delta$ (see Example
\ref{second_ex})
we define a 
particularly nice regular subdivision $\esd_r(T(\Delta)) =
\sd_\omega(rT(\Delta))$ of $rT(\Delta)$.  
The toric complex $\esd_r (T(\Delta))$ is simplicial and
the associated abstract simplicial complex is the {\em 
$r$-fold edgewise subdivision} $\esd_r(\Delta)$ of $\Delta$ defined in
\ref{edgewise_subdiv}. 
We show that the initial ideal 
$\ini_\omega(I_{rT(\Delta)})$ 
of the defining ideal $I_{rT(\Delta)}$ of the face ring of $rT(\Delta)$ is
generated 
by square-free monomials, and we specify a Gr\"obner basis for
$I_{rT(\Delta)}$ with respect to any monomial order $<$ refining the
weight order $\omega$.

Some ideas presented in this paper can be found elsewhere in
different contexts. First of all the concept of a toric complex is a variation
on the polyhedral complexes defined in \cite{Mumford}, and regular
subdivision appears there under the name projective
subdivision. 
There is a natural toric complex associated to simplicial fan, and
this toric complex plays an important r\^ole in the work 
of Cox on the homogeneous coordinate ring of a
toric variety \cite{Cox_Coordinate}. The work of Cox has 
motivated the concept of stacky 
fans, corresponding to a kind of toric complexes consisting of subsets
$F$ of a finitely generated abelian group instead of subsets of $\Z^d$,
considered in \cite{Borisov_Chen_Smith}.
The idea of considering arbitrary 
finite subsets $F$ of $\Z^d$ was promoted by Sturmfels and his coauthors, for
example in \cite{HOMAST} and \cite{ST}. They occasionally call $F$ a vector
configuration. 
In \cite{ST_Tohoku} and \cite{ST} Sturmfels identified the
radical of the initial ideal $\ini_\omega (\Ic_F)$ of $\Ic_F$
with the defining ideal of a Stanley-Reisner 
ring in the case where $\omega$ defines a triangulation of
$\cone(F)$ using Gr\"obner bases and integer programming methods. 
We make a different approach using direct geometric arguments that allows
us to extend Sturmfels' theorem to subdivisions of toric complexes and
face rings over    
commutative rings instead of over fields. 
The results of Section
\ref{thefacering} where we
examine the face ring of a toric complex are motivated by 
the paper \cite{EIST} of Eisenbud and Sturmfels.  
Finally, edgewise subdivision has been studied by many people. 
The basic idea is to subdivide a triangle along its edges:\\
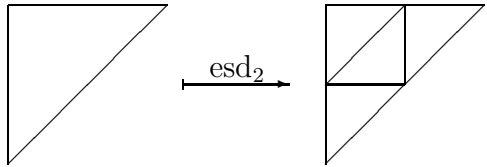
\begin{center}
\begin{picture}(180,60)
\put(0,0){
  \begin{picture}(60,60)
      \put(0,60){\line(1,0){60}}
      \put(0,0){\line(0,1){60}}
      \put(0,0){\line(1,1){60}}    
  \end{picture}
}
\put(70,28){\line(0,1){4}}
\put(70,30){\vector(1,0){40}}
\put(80,33){\shortstack{$\esd_2$}}
\put(120,0){
\begin{picture}(30,30)
  \put(0,30){\line(1,0){30}}
  \put(0,0){\line(0,1){30}}
  \put(0,0){\line(1,1){30}}
\end{picture}
}
\put(120,30){
\begin{picture}(30,30)
  \put(0,30){\line(1,0){30}}
  \put(0,0){\line(0,1){30}}
  \put(0,0){\line(1,1){30}}
\end{picture}
}
\put(150,30){
\begin{picture}(30,30)
  \put(0,30){\line(1,0){30}}
  \put(0,0){\line(0,1){30}}
  \put(0,0){\line(1,1){30}}
\end{picture}
}
\end{picture}
\\  
\end{center}
As observed by Freudenthal in \cite{Freudenthal}, as opposed to
barycentric subdivision, the pieces of
iterated edgewise subdivision do not become long and thin.
Knudsen and Mumford studied variations of edgewise subdivision of
polyhedral complexes in \cite{Mumford}. The name edgewise subdivision
was introduced by Grayson in \cite{Gr} where edgewise subdivision of
simplicial sets is used to obtain operations in higher algebraic $K$-theory.
Edgewise subdivision of cyclic sets is essential for the construction
of topological cyclic homology
\cite{BHM}.

The paper is organized as follows: in Section \ref{toriccomplexes}
we give our definition of a toric complex together with some examples, and we
characterize elementary properties of toric complexes.
In Section \ref{thefacering} we define the face ring of a toric complex, and
in the case where $K$ is a field we characterize the
$K$-algebras occurring as face-rings. 
Furthermore we
show that the face ring is compatible with gluing of toric complexes and we
give a construction on toric complexes corresponding to the Veronese subring
of a graded ring. 
In Section \ref{subdivisions} we introduce regular subdivisions of 
\toric complexes and we prove our main 
result Theorem \ref{main_result}.
In Section \ref{edgewisesubdivision} we discuss a particular regular
subdivision of the toric complex associated to an abstract simplicial
complex, which we 
call the edgewise subdivision. 

The authors are grateful to Prof.\ W.\ Bruns  
for inspiring discussions on the subject of the paper.

\section{Prerequisites}
In this section we fix some notation and recall some standard definitions.
Let $F$ be a finite subset of $\R^d$.
A convex combination of elements of $F$ 
is a sum $\sum_{a \in F} r_a a$ with  $0 \leq r_a$ for $a \in F$
and $\sum_{a \in F} r_a=1$.
The set of convex combinations 
$\conv(F)$ of elements of $F$ is called 
the {\em convex hull} of $F$. 
Similarly, a positive linear combination of elements of $F$
is a sum $\sum_{a \in F}^n r_a a$ with  $r_a \in \R_+ = \{x \in \R \colon x
\ge 0 \}$ for $a \in F$.
The set $\cone(F)$ of positive linear combinations 
of elements of $F$ is called 
the {\em cone} generated by $F$. 
By convention $\cone(\emptyset)=\{0\}$ and $\conv(\emptyset) = \emptyset$.

To a linear form $\alpha$ on  $\R^d$ and $c \in \R$ we 
associate the affine hyperplane 
$H_{\alpha}(c) = \alpha^{-1}(c)$
and the half-space 
$H_{\alpha}^-(c) = \alpha^{-1}((-\infty,c])$.
An intersection 
$P=\bigcap_{i=1}^n  H_{\alpha_i}^-(c_i)$
of finitely many half-spaces is called a {\em polyhedron}.
A {\em face} of a polyhedron $P$ is the intersection
of $P$ with an additional hyperplane $H_{\beta}(d)$ with the property that
$P \subseteq H_{\beta}^-(d)$.
A {\em polytope} is a bounded polyhedron. 
The set $\conv(F)$ is a polytope for every finite subset $F$ of $\R^d$ and
every polytope is of this form.
The polytope $\conv(F)$ is a {\em simplex} 
if every element in $\conv(F)$ has a unique representation
as a convex combination of the elements of $F$. 
A {\em cone} is a finite intersection of half-spaces of the form 
$H_{\alpha_i}^-(0)$.
The set $\cone(F)$ is a cone for every finite subset $F$ of $\R^d$ and every
cone is of this form. 
For the theory of polyhedrons
and related things we refer to the books of Schrijver \cite{SCH}
and Ziegler \cite{ZI}.

Given a commutative monoid $M$, the monoid algebra $K[M]$ is the set of
functions $f \colon M \to K$ with finite support with $(f+g)(m) =
f(m) + g(m)$ and $(fg)(m) = \sum_{m_1 + m_2 = m} f(m_1)g(m_2)$ for $m \in M$. 
Given a finite set $F$, 
the free commutative monoid on $F$ is the monoid $\N^F$ 
consisting of functions $u \colon F \to \N$, and with $(u+v)(a) = u(a)
+ v(a)$ for $u,v \in \N^F$ and $a \in F$. The polynomial ring $K[F]$ is
the monoid algebra $K[\N^F]$. A monomial in $K[F]$ is a function of the form
$x^u 
\colon \N^F \to K$ with $x^u(u) = 1$ and $x^u(v) = 0$ for $v \ne u$, and a
polynomial $f \in K[F]$ can be represented as the sum $f = \sum_{u
\in \N^F} f_u x^u$ where $f_u = f(u) \in K$.
For $a \in F$ we let $u_a \colon F \to \N$ denote the 
function defined by $u_a(a) = 1$ and $u_a(b) = 0$ for $b \ne a$ and we let
$x_a = x^{u_a}$.
If $G \subseteq F$, then we let $x_G$ be
the {\em square-free} monomial $\prod_{a \in G} x_a$.

Given a function $f \colon F \to G$    
of finite sets, the homomorphism
$f_* \colon K[F] \to K[G]$ is defined by $f_*(x_a) = x_{f(a)}$ for $x_a \in 
K[F]$, and the homomorphism
$f^* \colon K[G] \to K[F]$ is defined by
$f^*(x_b) = \sum_{f(a) = b} x_a$ for $x_b \in
K[G]$. Here a sum indexed on the empty set is
zero by convention.
The $r$-th Veronese subring 
of 
a $\Z$-graded ring $R = \oplus_{i \in \Z} R_i$ is the $\Z$-graded ring
$R^{(r)} = 
\oplus_{i \in \Z} R_{ri}$. 

\section{Toric complexes} 
\label{toriccomplexes}
In this section we introduce \toric complexes,
the objects studied in this paper.
\begin{defn}
\label{toric_complex}
  A {\em \toric complex} $T$ consists of:
  \begin{enumerate}
    \item[(a)] 
    A finite partially ordered set $(\Pi_T,\subseteq)$ consisting of finite
    sets, ordered by inclusion,
    \item[(b)] 
    for every $F \in \Pi_T$ an injective function 
    $T_F\colon F \to \Z^{d_F}\setminus \{0\}$ for some $d_F > 0$, 
    \item[(c)] for every pair $F, G \in \Pi_T$ with $F \subseteq G$ an
    injection $T_{FG} \colon \Z^{d_F} \to \Z^{d_G}$ of abelian groups, 
  \end{enumerate}
  subject to the following conditions:
  \begin{enumerate}
    \item if $F,G \in \Pi_T$, then $F \cap G \in \Pi_T$,
    \item for all $F, G \in \Pi_T$ with $F \subseteq G$ and $a
    \in F$ we have $T_{FG}(T_F(a)) = T_G(a)$, 
    \item for every triple $F, G, H \in \Pi_T$ with $F \subseteq G
    \subseteq H$ we have $T_{GH} \circ T_{FG} = T_{FH}$,
    \item if $G \in \Pi_T$ and $F \subseteq G$, then $F \in \Pi_T$ if and
    only if $\cone(T_G(F))$ is a face of $\cone(T_G(G))$ satisfying
    $T_G(F) = \cone(T_G(F)) \cap T_G(G)$.
\end{enumerate}
The {\em set of generators} of $T$ is the union $\Gen(T) = \cup_{F\in \Pi_T} F$
and 
the {\em faces} of $T$ are the elements $F \in \Pi_T$. If $F$ and
$G$ are faces of $T$ with $F \subseteq G$ we say that $F$ is a {\em face}
of $G$ in $T$.
For a face $F\in \Pi_T$ we define 
the monoid homomorphism $T_F\colon \N^F \to \Z^{d_F}$,
$T_F(u) =\sum_{a \in F} u(a)T_F(a)$. 
\end{defn}
A toric complex $T$ is {\em embedded} if
there exist $d \in \N$ such that 
$T_{FG}$ is the identity on $\R^d$ for every $F,G \in \Pi_T$ with 
$F\subseteq G$.
In particular, we have $T_F(a)=T_G(a)$ for $a\in F$.
Note that for a \toric complex $T$ the collection $\{ \cone(F) \colon
F \in \Pi_T\}$ is 
a polyhedral complex as considered for example in \cite{Mumford} 
consisting of rational cones.
If $T$ is embedded, then this is a fan in the lattice
$\Z^d$.

Whenever appropriate we shall consider a partially ordered set
$\Pi$ as a category with objects set equal to the underlying set of
$\Pi$ and with morphism sets $\Pi(F,G)$ consisting of exactly one element
$F \to G$ if $F \le G$ and with $\Pi(F,G)$ empty otherwise. 
We let $\Pi^{\op}$ denote the opposite category of $\Pi$ with the same
set of objects as $\Pi$ and with $\Pi^{\op}(G,F) = \Pi(F,G)$.

Given a toric complex $T$ and elements $F \subseteq G$ of $\Pi_T$,
there are continuous injective maps 
$\conv(T_F(F)) \to \conv(T_G(G))$ and
$\cone(T_F(F)) \to \cone(T_G(G))$
induced by the injective homomorphism $T_{FG}$. Thus we can consider
$F \mapsto \conv(T_F(F))$ and $F \mapsto \cone(T_F(F))$ as functors from
 $\Pi_T$ 
to the category of topological
spaces.
We define the spaces $|T|$ and $||T||$ as the colimits:
$$
|T| = \colim{F \in \Pi_T} \conv(T_F(F)) \text{ and }
||T|| = \colim{F \in \Pi_T} \cone(T_F(F)).
$$  
Observe that if $T$ is embedded, then the topological spaces 
$|T| = \bigcup_{F \in \Pi_T} \conv(T_F(F))$ 
and  
$||T|| = \bigcup_{F \in \Pi_T} \cone(T_F(F))$
are unions of subspaces of $\R^d$.
\begin{ex}
\label{first_ex}
A subset $G$ of $\Z^d$ gives
rise to an embedded \toric complex $T(G)$. 
The function $T(G)_G$ is the
inclusion of $G$ in $\Z^d$ and $\Pi_{T(G)}$ consist of those subsets
$F$ of $G$ satisfying axiom (iv) of Definition \ref{toric_complex}.

This type of \toric complexes correspond to vector configurations
as studied for example in \cite{HOMAST} or \cite{OSRT}.
\end{ex}
\begin{ex}
\label{second_ex}
Let $\Delta$ be an abstract simplicial complex on the vertex set 
$V=\{1,\dots,d+1\}$, i.e.\
$\Delta$ is a set of subsets of $V$ 
and $F\subseteq G \in \Delta$ implies $F \in \Delta$.
Let $e_1,\dots,e_{d+1}$ denote the standard generators of $\Z^{d+1}$
and define $\iota_V \colon V \to \Z^{d+1}$ by $\iota_V(i) = e_i + e_{i+1} +
\dots + e_{d+1}$ for $1 \le i \le d+1$.
We associate an embedded 
\toric complex $T(\Delta)$ to $\Delta$ as follows: 
\begin{enumerate}
\item[(a)] 
$\Pi_{T(\Delta)}$ is the partially ordered set $\Delta$,
\item[(b)] 
for $F \in \Delta$ we let $T(\Delta)_F$ denote the restriction of
$\iota_V$ to $F \subseteq V$,
\end{enumerate}
Observe that $|T(\Delta)|$ is homeomorphic to
the usual geometric realization of $\Delta$.
\end{ex}
\begin{defn}
\label{define_properties}
Let $T$ be a toric complex.
\begin{enumerate}
\item
$T$ is {\em $\Z$-graded} if $T_{FG}\colon \Z^{d_F} \to \Z^{d_G}$ preserves
 the last coordinate for all $F, G \in \Pi_T$ with $F \subseteq G$.  
\item
$T$ is {\em $\N$-graded} if it is $\Z$-graded and the last coordinate of
$T_F(a)$ is positive for all $a \in F$ and $F \in \Pi_T$.
\item 
$T$ is {\em standard graded} if it is $\Z$-graded and the last
coordinate of 
$T_F(a)$ is $1$ for all $a \in F$ and $F \in \Pi_T$.
\item
$T$ is {\em pointed} if $\emptyset \in \Pi_F$.
\item  
$T$ is {\em simplicial} if $\Pi_T$ is an abstract simplicial complex.
\end{enumerate}
\end{defn}
Observe that the \toric complex $T(\Delta)$ of Example \ref{second_ex} is simplicial.
\begin{ex}
  \label{fanex}
  Let $\Sigma$ be a rational fan in a lattice $N$, that is, $\Sigma$
  is a finite collection of rational polyhedral cones
  $\sigma$ in $N_{\R} = N \otimes_{\Z} \R$, satisfying: every face of
  a cone in $\Sigma$ is also a cone in $\Sigma$, and the intersection
  of two cones in $\Sigma$ is a face of each. 
  Choosing a set
  $G_{\sigma} \subseteq N\setminus \{0\}$ of generators for each cone
  $\sigma$ in $\Sigma$ and an isomorphism $\varphi \colon N \to \Z^d$
  we can construct an embedded toric complex
  $T(\Sigma,(G_{\sigma})_{\sigma \in \Sigma}, \varphi)$ as follows:
  for $\sigma 
  \in \Sigma$ we let $F_{\sigma}$ denote the union of the sets
  $\varphi(G_{\gamma})$ of images of generators of faces $\gamma$ of $\sigma$
  under $\varphi$ and we
  define $\Pi_{T(\Sigma,(G_{\sigma})_{\sigma \in \Sigma}, \varphi)}$
  to be the collection of subsets of $\Z^d$ of the form $F_{\sigma}$ for
  $\sigma 
  \in \Sigma$. 

  There is a preferred toric complex associated to a simplicial fan,
  namely the one where we take $G_\rho$ to be the set of unique
  generators $\rho \cap N$ if $\rho$ is a ray in $\Sigma$, and where
  $G_\sigma$ is the union over the rays $\rho$ in $\sigma$ of the sets
  $G_\rho$. On the other hand, the underlying fan $\{\cone(T_F(F)) \colon F
  \in \Pi_T\}$ of an embedded toric complex $T$ is simplicial if $T$ is
  simplicial.  If we want to subdivide 
  $T(\Sigma,(G_{\sigma})_{\sigma \in \Sigma}, \varphi)$, it is
  necessary to choose the sets $G_{\sigma}$ different from the
  preferred ones described above. In Definition
  \ref{define_multiple} below we present a way to do this. 
\end{ex}
\begin{ex}
Lattice polyhedral complexes in the sense of 
\cite[Definition  2.1]{BRGU} 
correspond to the subclass of 
the class of standard graded toric complexes
consisting of toric complexes $T$ with
the 
property that $F$ is the set of vertices of the polytope
$\conv(T_F(F))$
for every $F \in \Pi_T$. For a concrete
example of a non-embeddable toric complex we refer to
\cite[Proposition 2.3]{BRGU}.
\end{ex}

\begin{lem}
\label{helper1}
Let $T$ be a \toric complex. 
\begin{enumerate}
\item
$T$ is pointed if and only if zero is a vertex of
$\cone(T_F(F))$ for every $F \in \Pi_T$.
\item
If $T$ is $\N$-graded, then it is pointed. 
\item
$T$ is simplicial if and only
if for all $F\in \Pi_T$ the elements $T_F(a)$ for
$a \in F$ are linearly independent.
\end{enumerate} 
\end{lem}
\begin{proof}
We only prove 
(iii).
Assume that
for all $G \in \Pi_T$ the elements $T_G(a)$ for
$a \in G$ are linearly independent.
If $F \subseteq G$ and $G\in \Pi_T$,
then 
it is easy to see that
$\cone(T_G(F))$ is a face of $\cone(T_G(G))$
and $T_G(F) = T_G(G) \cap \cone(T_G(F))$.
Since $T$ is a \toric complex this implies $F \in \Pi_T$
and thus $\Pi_T$ is a simplicial complex.

On the other hand, suppose that $\Pi_T$ is an abstract simplicial complex
and let $G \in \Pi_T$.
Assume that the elements of $T_G(G)$ are linearly dependent.
Thus for some $a \in G$  there is a relation
$T_G(a)=\sum_{b\in G, b \neq a} \lambda_b T_G(b)$ with $\lambda_b \in \R$.
By our assumption $F=G\setminus \{a\}$ generates a proper face
$\cone(T_G(F))$ of $\cone(T_G(G))$
with $T_G(F)=\cone(T_G(F)) \cap T_G(G)$.
Choosing a linear form $\alpha$
on $\R^{d_G}$ with $\cone(T_G(F)) = \alpha^{-1}(0) \cap \cone(T_G(G))$
we get the  
contradiction 
$0\neq \alpha(T_G(a))=\sum_{b\in F} \lambda_b \alpha(T_G(b))=0.$
Hence the elements in  $T_G(G)$ are linearly independent.
\end{proof}
\begin{defn}
We call $S$ a {\em subcomplex} of a toric complex $T$ and write
$S \subseteq T$ 
if $S$ is a \toric complex 
with $\Pi_S \subseteq \Pi_T$ such that
$S_F = T_F$ and $S_{FG} = T_{FG}$ 
for $F,G\in \Pi_S$ with $F \subseteq G$. 
\end{defn}
Observe that by the definition of a \toric complex
we have for $G \in \Pi_S$ and $F \subseteq G$
that $F \in \Pi_S$ if and only if $F \in \Pi_T$.

\section{The face ring of a toric complex}
\label{thefacering}
Stanley associated in \cite{ST87} a $K$-algebra to a rational fan.
In this section we translate Stanley's definition to the situation of \toric
complexes.  
Recall the notation introduced in Section 1 and 2 
which will be used in the following.
\begin{defn}
\label{face_ring}
  Let $T$ be a toric complex. The
  {\em face ring} $K[T]$ is the $K$-algebra
  $K[T] = K[\Gen(T)]/I_T$,
  where $I_T =  J_T + \sum_{F \in \Pi_T} I_F$. Here $J_T$ is the monomial
  ideal such that a monomial $x^u$ is {\em not} in $J_T$ if and only if there
  exists $F \in \Pi_T$ with 
  $\supp (x^u) \subseteq F$. For $F \in \Pi_T$
  the ideal $I_F$ is generated by the image of the kernel $\Ic_F \subseteq
  K[F]$ of the 
  homomorphism  
  $K[F] \to K[\Z^{d_F}]$ induced by $T_F$
  under the inclusion $K[F] \to K[\Gen(T)]$.
\end{defn}

Recall that an ideal in a polynomial ring over $K$ is {\em binomial},
if it is generated by binomials $x^u-x^v$ and monomials $x^u$.
Such binomial ideals are strict in the sense that we do not
allow the  
generators to have coefficients different from $1$.
A {\em monomial} ideal is an ideal generated by monomials.

Observe that the ideal $I_T$ of a \toric complex is binomial, 
because for every $F \in \Pi_T$ the
ideal $I_F$ is generated by the binomials $x^u-x^v$ where 
$u,v \colon F \to \N$ satisfy that $T_F(u) = T_F(v)$
(e.g. see \cite[Lemma 4.1]{ST}).
Note that if $I_F$ is a monomial ideal, then $I_F = 0$.
If $T$ is $\Z$-graded, then $K[T]$ inherits a $\Z$-grading, and
$K[T]$ is $\N$-graded if $T$ is $\N$-graded. If $T$ is standard graded,
then $K[T]$ is generated by elements of degree one.

\begin{thm}
\label{toric_proj}
Let $S$ be a subcomplex of a \toric complex $T$ and
let $i \colon \Gen(S) \to \Gen(T)$ 
denote the inclusion of the generators of $S$ in the generators of $T$.
\begin{enumerate}
\item
The homomorphism $i^* \colon K[\Gen(T)] \to K[\Gen(S)]$ induces a  
surjective homomorphism
$
i_{ST}^* \colon K[T] \to K[S]
$
with  
$$
\Ker(i_{ST}^*)=
(x_F \colon F\subseteq \Gen(T) \text{ is not contained in any } G \in \Pi_S).
$$
\item
The homomorphism $i_* \colon K[\Gen(S)] \to K[\Gen(T)]$ induces an injective
homomorphism ${i_{ST}}_* \colon K[S] \to K[T]$ with 
$i_{ST}^* \circ {i_{ST}}_* = \id_{K[S]}$. 
\end{enumerate}
\end{thm}

\begin{proof}
To prove that $i_{ST}^*$ is well defined
we need to check that $i^*$ maps the ideal 
$I_T = J_T + \sum_{F \in \Pi_T} I_F$ into the ideal $I_S = J_S + \sum_{G \in
\Pi_S} I_G$. 

Clearly $J_T$ is mapped into $J_S$.
For $F \in \Pi_T$ the ideal $I_F$ is
generated by binomials $x^u - x^v$ where $x^u$ and $x^v$ have
support contained in $F$ and satisfy that $T_F(u) = T_F(v)$.
For $G\subseteq F$ such that $\cone(T_F(G))$ is a face of $\cone(T_F(F))$
we have $\supp(x^u) \subseteq G$ if and only if $T_F(u) \in
\cone(T_F(G))$. Thus $\supp(x^u) \subseteq G$ if and only if
$\supp(x^v) \subseteq G$. 
The element $i^*(x^u)$ is {\em not} in $J_S$ if and only if there exists a
$G \in \Pi_S$ such that $\supp(x^u) \subseteq G$. 
Hence $i^*(x^u)$ is {\em not} in $J_S$ if
and only if the same is the case for $i^*(x^v)$. 
Now if $i^*(x^u)$ and $i^*(x^v)$ are in $J_S$ we are
done. Otherwise there exists $G \in \Pi_S$ such that $\supp(x^u)$ and
$\supp(x^v)$ are contained in $G$ 
and therefore 
$i^*(x^u-x^v)$ is in $I_G$.
We conclude that $i^*(I_T) \subseteq I_S$
and we get the induced homomorphism $i_{ST}^*\colon K[T]\to K[S]$.
This homomorphism is clearly surjective and has the described kernel.

Analogously one shows that 
${i_{ST}}_* \colon K[S] \to K[T]$ is an injective homomorphism
with $i_{ST}^* \circ {i_{ST}}_* = \id_{K[S]}$. This concludes the proof.
\end{proof}

%
%

We need the following variation of Corollary 2.4 in \cite{EIST}.
\begin{lem}
\label{binomial_helper}
Let $K$ be a field and let $P$ be a prime ideal
in $K[F]$ generated by binomials $x^u - x^v$ for some finite set
$F$. There exists a unique direct summand $L_P$ of $\Z^F$ such that the
composition $K[F]= K[\N^F] \subseteq K[\Z^F] \to K[\Z^F/L_P]$ induces 
an embedding $K[F]/P \subseteq K[\Z^F/L_P]$.  
\end{lem}
\begin{proof}
Since
$x^{m+n} - 1 = x^m(x^{n}- 1) + (x^m -1)$ the elements $m \in \Z^F$
satisfying $x^m - 1 \in PK[\Z^F]$ form a lattice $L_P$ with 
$PK[\Z^F]$ contained in the kernel $I(L_P)$ 
of the homomorphism $\gamma\colon K[\Z^F] \to K[\Z^F/L_P]$. 
Since
$P=PK[\Z^F]\cap K[\N^F]$,
it suffices to prove that the
homomorphism $\delta\colon K[\Z^F]/PK[\Z^F] \to K[\Z^F/L_P]$
is an isomorphism and that
$L_P$ is a direct summand of $\Z^F$.
If $\delta$ is an isomorphism, 
then $L_P$ must be a direct summand of $\Z^F$
because $PK[\Z^F]$ is a prime ideal.

An element
$
f 
= \sum_{u \in \Z^F} f_u x^u 
$
is mapped to 
$$
\gamma(f)=\sum_{v + L_P  \in \Z^F/L_P}  (\sum_{m \in L_P} f_{v+m}) x^{v+L_P}.
$$
We have that $\gamma(f)=0$ if and only if $\sum_{m \in L_P} f_{v+m}=0$  
for all $v+L_P \in \Z^F/L_P$  and this implies that 
$\sum_{m \in L_P} f_{v+m} x^{m}$ is in the ideal generated by the $x^m-1$.
Thus we see that the ideal
$I(L_P)$ is generated by the elements $x^m - 1$ for $m \in L_P$.
We conclude that $I(L_P) \subseteq PK[\Z^F]$,
and that the homomorphism $\delta$ in question is an isomorphism.
This concludes the proof.
\end{proof}

Let $\Inew\subset K[x_1,\dots,x_n]$ 
be a  binomial ideal.
For a subset 
$F \subseteq \{1,\dots,n\}$ 
let $P_F= \Inew \cap K[F]$.
As in \cite[Cor. 1.3]{EIST} we see that $P_F$ is a  binomial ideal in $K[F]$.
Let $\Pi_{T(\Inew)}$ denote the collection of subsets $F$ of $\{1,\dots,n\}$
such that 
$P_F$ is a prime
ideal not containing any monomial 
and such that the projection $K[x_1,\dots,x_n] \to K[F]$
induces a homomorphism $K[x_1,\dots,x_n]/I \to K[F]/P_F$.
\begin{thm}
\label{char_face}
Let $K$ be a field, $X = \{1,\dots,n\}$ and
$\Inew \subset K[X] = K[x_1,\dots,x_n]$ be a  binomial ideal.
There exists a \toric complex $T(\Inew)$ with face set $\Pi_{T(\Inew)}$ 
and $K[T(\Inew)] \cong K[X]/\bigcap_{F \in \Pi_{T(\Inew)}} (P_F + (x_i \colon i \notin F))$. 
In particular, 
$\Inew=\bigcap_{F \in \Pi_{T(\Inew)}} (P_F + (x_i \colon i \notin F))$ if and
only if the natural 
homomorphism $K[X]/\Inew \to K[T(\Inew)]$ is an isomorphism.
\end{thm}
\begin{proof}
Let us first construct an abstract toric
complex $T'(\Inew)$ in the sense that the targets for the maps $T'_F$ for $F
\in \Pi_{T'(\Inew)} = \Pi_{T(\Inew)}$ are abstract finitely generated free
abelian groups 
instead of abelian groups of the form $\Z^{d_F}$.
For every $F \in \Pi_{T(I)}$ there is an isomorphism $K[X]/\Inew + (x_i \colon i
\notin F) \cong K[F]/P_F$. Let $L_{P_F}$ be as in Lemma \ref{binomial_helper}
and
let $T'_F$ denote the composition 
$F \subseteq \N^F \subseteq \Z^F \to \Z^F/L_{P_F}$.
If $F \subseteq G$, then the inclusion 
$\Z^F \subseteq \Z^G$
and the fact $P_F =P_G \cap K[F]$ 
induce an inclusion $T'_{FG} \colon
\Z^F/L_{P_F} \to \Z^{G}/L_{P_G}$
.
Choosing isomorphisms
$\Z^F/L_{P_F} \cong \Z^{d_F}$ and
applying \cite[Theorem 6.1.7]{BRHE} it is straight forward to check
that we obtain a toric complex $T(\Inew)$. 

In order to identify $K[T(I)]$ we assume without loss of generality
that $\Gen(T(I)) = X$.
Note that for $F \in \Pi_{T(\Inew)}$ the ideal $P_F$ agrees with the kernel
$\Ic_F$ of the
homomorphism $K[F] \to K[\Z^{d_F}]$ induced by $T_F$. 
Considering the inclusions of \toric complexes
$T(F) \subseteq T(I)$, Theorem \ref{toric_proj} ensures that
$i_{*}i^*(f) \in I_F$ for $f \in I_{T(I)}$.
Since $f-i_{*}i^*(f) \in (x_i \colon i \notin F)$ and
$f$ was arbitrary,
we get that $I_{T(I)}$ 
is contained in 
$(P_F + (x_i \colon i \notin F))$.
Thus
$I_{T(\Inew)}$ is contained in
$\bigcap_{F \in \Pi_{T(\Inew)}} 
(P_F+ (x_i \colon i \notin F))
$.

Let $f\in \bigcap_{F \in \Pi_{T(\Inew)}} (P_F + (x_i \colon i \notin F))$ 
and  let $F_1,\dots,F_m$ be the maximal faces of $T(\Inew)$.
Recall that the inclusion $i$ of $F_1$ in  $X$ induces
homomorphisms $i^* \colon K[X] \to K[F_1]$ and 
$i_* \colon K[F_1] \to  K[X]$.
The element
$f - i_* i^*(f) \in K[X]$ is in the ideal
$\left( \bigcap_{F \in \Pi_{T(\Inew)}} (P_F + (x_i
\colon i 
\notin F))\right) \cap (x_i \colon i \notin F_1)$ 
because $i_*i^*(f) \in I_F \subseteq I_{T(\Inew)}$. 
Proceeding by induction on $m$ we 
find an element $g$ of $I_{T(\Inew)}$ such that 
$f - g \in 
\bigcap_{1 \le j \le m} (x_i \colon i \notin F_j) = J_{T(\Inew)} 
\subseteq I_{T(\Inew)}$.
Hence $f \in I_{T(\Inew)} 
$ and it follows that 
$I_{T(\Inew)}
=\bigcap_{F \in \Pi_{T(\Inew)}} (P_F + (x_i \colon i \notin F))$
and that 
$K[T(\Inew)]=K[X]/\bigcap_{F \in \Pi_{T(\Inew)}} (P_F + (x_i
\colon i \notin F)).$
The last claim of the theorem is a consequence of this fact.
\end{proof}
If a \toric complex $T$ is simplicial, then all the ideals $I_F=0$
for $F \in \Pi_T$ 
because 
the linear independence of the elements of $T_F(F)$ implies that
the ring $K[M_{T_F(F)}]$ is a polynomial ring.
Hence $I_T=J_T$ is a square-free monomial ideal 
and the face ring $K[T]$ 
coincides with the so-called Stanley-Reisner ring of the abstract
simplicial complex $\Pi_T$.
(See for example \cite{Stanley_book} for more details on this subject.)
We show that the converse is also true.
\begin{prop}
\label{helper3}
A toric complex $T$ is simplicial if and only if 
$I_T$ is a monomial ideal.
In this case
$K[T]$ is the Stanley-Reisner ring 
of the abstract simplicial complex $\Pi_T$.
\end{prop}
\begin{proof}
If a \toric complex $T$ is simplicial, then $I_T$ is a monomial ideal
as shown above. 

Therefore assume that $I_T$ is a monomial ideal. 
Given $F \in \Pi_T$ we consider the toric subcomplex $T(F)$ of Example
\ref{first_ex}. Let $i\colon F=\Gen(T(F)) \to \Gen(T)$ be the inclusion. 
Since $I_T$ is monomial, so is the ideal $I_{T(F)} = i^* (I_T)$.
It follows that $I_{T(F)}=0$, 
because $K[F]/I_{T(F)}$ is an affine monoid ring.
Thus
the elements of $T_F(F)$ are linearly independent. 
By \ref{helper1} (iii) we conclude that $T$ is simplicial.
Then $I_T=J_T$ and $K[T]$ is the Stanley-Reisner ring of $\Pi_T$.   
\end{proof}

Given a \toric complex $T$ and subcomplexes $Q,S \subset T$, there is a
subcomplex $R=Q \cap S$ of $T$ defined by letting 
$\Pi_{R} =\Pi_Q \cap \Pi_S$. 
One can show that the space $|R|=|Q \cap S|$ is
isomorphic to $|Q| \cap |S|$.
We write 
$T = Q \cup S$ if 
$\Pi_{T} = \Pi_Q \cup \Pi_S$. 
Since there are homomorphisms 
$i_{RQ}^* \colon K[Q] \to K[R]$ 
and $i_{RS}^* \colon K[S] \to K[R]$ we can consider the fiber product
$K[Q]\times_{K[R]} K[S]$.

Observe that the homomorphisms ${i_{ST}}_*$ and $i_{ST}^*$ are {\em twins} in the
sense of Notbohm and Ray \cite{notbohm_ray}, 
that is
$
{i_{QT}}^* {i_{ST}}_* = {i_{RQ}}_* {i_{RS}}^*
\text{ and }
{i_{ST}}^* {i_{QT}}_* = {i_{RS}}_* {i_{RQ}}^*.
$
The face ring of $T$ and the face rings of $Q$ and $S$
are related as follows:
\begin{prop}
\label{fibre}
If $T = Q \cup S$, then the homomorphism 
$$
(i_{QT}^*,i_{ST}^*) \colon K[T] \to K[Q]\times_{K[Q \cap S]} K[S]
$$ 
is an isomorphism.
\end{prop}
\begin{proof}
  As above let $R = S
  \cap Q$.
  We prove that the
  additive homomorphism 
  $\beta \colon K[Q]\times_{K[R]} K[S] \to K[T]$ defined by
  $
  \beta(a,b) = {i_{QT}}_*(a) + {i_{ST}}_*(b) - {i_{RT}}_*( {i_{RQ}}^*(a)) 
  $
  for $(a,b) \in K[Q]\times_{K[R]} K[S]$
  is a homomorphism of
  rings. 
  Using the fact that 
  ${i_{ST}}_*$ and $i_{ST}^*$ are twins
  we see that $\beta$ is an additive inverse to
  $(i_{QT}^*,i_{ST}^*)$.

  The only thing left to check is that
  $\beta(a,b) \beta(a',b') = \beta(aa',bb')$. 
  Since the support $\supp(x^u)$ of a monomial $x^u$ in the polynomial
  $({i_{QT}}_*(a) - {{i_{RT}}_*}{i_{RQ}^*}(a))$ satisfies that 
  $\supp(x^u) \cap (\Gen(Q) \setminus \Gen(R)) \ne \emptyset$, 
  and since 
  $\supp(x^v) \cap (\Gen(S) \setminus \Gen(R)) \ne \emptyset$ 
  for a monomial $x^v$ in the polynomial 
  $({{i_{ST}}_*}(b') - {{i_{RT}}_*}{i_{RS}^*}(b'))$
  we have that  
$$
({i_{QT}}_*(a) - {{i_{RT}}_*}{i_{RQ}^*}(a))
({{i_{ST}}_*}(b') - {{i_{RT}}_*}{i_{RS}^*}(b')) \in J_T.
$$ 
Similarly we see that 
$
({{i_{ST}}_*}(b) - {{i_{RT}}_*}{i_{RS}^*}(b))
({i_{QT}}_*(a') - {{i_{RT}}_*}{i_{RQ}^*}(a'))
\in J_T,$ 
and we compute that
  \begin{eqnarray*}
    && {\beta(aa',bb')} \\
    &=& 
\beta(aa',bb') + 
({i_{QT}}_*(a) - {{i_{RT}}_*}{i_{RQ}^*}(a))
({{i_{ST}}_*}(b') - {{i_{RT}}_*}{i_{RS}^*}(b')) \\
&& 
\quad + 
({{i_{ST}}_*}(b) - {{i_{RT}}_*}{i_{RS}^*}(b)) 
({i_{QT}}_*(a') - {{i_{RT}}_*}{i_{RQ}^*}(a'))\\ 
    &=& {i_{QT}}_*(aa') + {{i_{ST}}_*}(bb')
     -{{i_{RT}}_*}{i_{RQ}^*}(aa') \\
    && 
\quad + 
{i_{QT}}_*(a){{i_{ST}}_*}(b') 
- 
{i_{QT}}_*(a) {{i_{RT}}_*}{i_{RS}^*}(b') - 
{{i_{RT}}_*}{i_{RQ}^*}(a) {{i_{ST}}_*}(b') \\
&& \quad + {{i_{RT}}_*}{i_{RQ}^*}(a) {{i_{RT}}_*}{i_{RS}^*}(b')\\
&& \quad + {{i_{ST}}_*}(b){i_{QT}}_*(a') -  
{{i_{RT}}_*}{i_{RS}^*}(b){i_{QT}}_*(a')  -  
{{i_{ST}}_*}(b){{i_{RT}}_*}{i_{RQ}^*}(a') \\
&& \quad + {{i_{RT}}_*}{i_{RS}^*}(b){{i_{RT}}_*}{i_{RQ}^*}(a') \\    
&=&        
{i_{QT}}_*(a){i_{QT}}_*(a') + 
{i_{QT}}_*(a){{i_{ST}}_*}(b') - 
{i_{QT}}_*(a){{i_{RT}}_*}{i_{RQ}^*}(a')\\
&& \quad + {{i_{ST}}_*}(b){i_{QT}}_*(a') + {{i_{ST}}_*}(b){{i_{ST}}_*}(b') - {{i_{ST}}_*}(b){{i_{RT}}_*}{i_{RQ}^*}(a') \\
&& \quad - {{i_{RT}}_*} {i_{RQ}^*}(a){i_{QT}}_*(a')- {{i_{RT}}_*}{i_{RQ}^*}(a) {{i_{ST}}_*}(b') + {{i_{RT}}_*}{i_{RQ}^*}(a){{i_{RT}}_*}{i_{RQ}^*}(a') \\    
&=& 
({i_{QT}}_*(a) + {{i_{ST}}_*}(b) - {{i_{RT}}_*}{i_{RQ}^*}(a))
({i_{QT}}_*(a') + {{i_{ST}}_*}(b') - {{i_{RT}}_*}{i_{RQ}^*}(a'))\\
&=& \beta(a,b) \beta(a',b').
\end{eqnarray*}
  This concludes the proof.
\end{proof}

Proposition \ref{fibre} can be generalized to the following situation.
\begin{thm}
\label{fibre_second}
Let $T$ be a \toric complex, $T_1,\dots,T_r$ be subcomplexes of
$T$ such that $T = T_1 \cup \dots \cup T_r$ and let $\Pc(r)$ denote
the partially ordered set consisting of all subsets of $\{1,\dots,r\}$
ordered by inclusion. Given $I \in \Pc(r)$ we let $T_I$ denote the
subcomplex $T_I = \cap_{i \in I} T_i$ of $T$. Then the natural homomorphism 
$$ 
K[T] \to \prolim{I \in \Pc(r)} K[T_I]
$$ 
is an isomorphism.
\end{thm}
\begin{proof}
  We prove the lemma by induction on $r$, the case $r = 1$ being
  obvious. Assume that the lemma holds for $r-1$. 
  In particular, the homomorphisms 
$$
K[T_1 \cup \dots \cup T_{r-1}] \to \prolim{I\in \Pc(r-1)} K[T_I]
\text{ and } 
K[(T_1 \cup \dots \cup T_{r-1})\cap T_r] \to \prolim{I\in \Pc(r), r \in I }K[T_I]
$$ 
are isomorphisms. 
Since $T = (T_1 \cup \dots \cup T_{r-1}) \cup T_r$ there is a chain of isomorphisms
  \begin{eqnarray*}
    K[T] &\cong& K[T_1 \cup \dots \cup T_{r-1}] \times_{K[(T_1 \cup
    \dots \cup T_{r-1}) \cap T_r]} K[T_r] \\
    &\cong& \prolim{I \in \Pc(r), r \notin I} K[T_I] 
    \times_{\prolim{I\in \Pc(r), r \in I}K[T_I]} K[T_r] \\
    &\cong& \prolim{I \in \Pc(r)} K[T_I].
  \end{eqnarray*}
  where the first is the isomorphism of \ref{fibre}, the
  second isomorphism is induced by the isomorphisms that hold by
  the inductive hypothesis and the third isomorphism uses the fact 
  that the objects in question have the same universal property. 
\end{proof}

For $F \in \Pi_T$ consider the complex $T(F)$ of Example \ref{first_ex}.
The face ring of $T(F)$ has the form $K[T(F)]=K[F]/\Ic_F$
where $\Ic_F$ is the kernel of the homomorphism $K[F] \to K[\Z^{d_F}]$. 
The following corollary was 
proved in \cite[Prop. 4.8]{EIST}
by considering the limit as a subset of a product.
\begin{cor}
\label{fibre_third}
Let $T$ be a \toric complex. Then
$$
K[T] 
\cong 
\prolim{F \in \Pi_T^{\op}} 
K[T(F)]
$$ 
\end{cor}
\begin{proof}
Note that $T(F) \cap T(G) = T(F \cap G)$ and apply \ref{fibre_second}
to the family 
$\{T(F)\}_{F \in \Pi_T}$. 
\end{proof}
Thus we have another representation of the face ring of a \toric complex.
For example this implies easily that $K[T]$ is reduced,
because $K[T(F)]$ is reduced for $F \in \Pi_T$
and the limit of reduced rings is reduced.
The next goal will be to interpret the Veronese subrings 
of the face ring
of a standard graded \toric complex.
The following construction will be used later to
subdivide simplicial toric complexes.
\begin{defn}
\label{define_multiple}
Given a \toric complex $T$ and $r \ge 1$ we let $rT$
consist of the data: 
\begin{enumerate}
\item[(a)] 
$\Pi_{rT}$ is the partially ordered set consisting of the
sets  
$$
rF=\{ u \in \N^F \colon \sum_{a \in F} u(a) = r \} / \sim
$$ 
where $F \in \Pi_T$ and $u \sim v$ if $T_F(u) = T_F(v)$.
\item[(b)]
For $rF \in \Pi_{rT}$ and $u \in rF$ let
$rT_{rF}(u) = T_F(u)$.
\item[(c)]
For $rF,rG \in \Pi_{rT}$ with $rF \subseteq rG$
let $rT_{rF,rG} = T_{F,G}$.
\end{enumerate}
\end{defn}
\begin{prop}
\label{verones_first}
Let $T$ be a \toric complex and $r \geq 1$ a positive integer.
Then $rT$ is a \toric complex.
\end{prop}
\begin{proof}
The sets $rF \in \Pi_{rT}$ are all subsets of the set 
$\{u \colon \Gen(T) \to \N \}/\sim$ 
where $u \sim v$ if there exists  
$G \in \Pi_T$ such that both $\supp(u)$ and $\supp(v)$ are contained in $G$
and $T_G(u) = T_G(v)$.
By construction $rT_{rF}\colon  rF \to \Z^{d_F} \setminus \{0\}$ is
injective and $rT_{rFrG} = T_{FG}$ is an injection of abelian groups. 
We have to verify conditions (i)-(iv) of the definition of a \toric complex.
Condition (i) holds because we have that 
$rF \cap rG = r(F \cap G)$ for $F, G \in \Pi_T$. Conditions (ii) and (iii) are
satisfied by construction, and condition (iv) holds, 
because $\cone(T_F(F)) = \cone(rT_{rF}(rF))$ for $F \in \Pi_T$. 
\end{proof}
Clearly if $T$ is a \toric complex 
and $S$ is a subcomplex of $T$ then
$rS$ is a subcomplex of $rT$.
\begin{thm}
\label{toric_veronese}
Let $T$ be a standard graded \toric complex and $r \geq 1$ a positive
number. Then $K[rT]=K[T]^{(r)}$ is the $r$-th Veronese ring of $K[T]$.
\end{thm}
\begin{proof}
We compute
$$
K[rT] 
\cong    
\prolim{rF \in \Pi_{rT}^{\op} } K[(rT)(rF)] 
\cong \prolim{F \in \Pi_T^{\op} } K[T(F)]^{(r)} 
\cong K[T]^{(r)}
$$
where the second isomorphism follows from the definition
of the \toric complex $rT$.
\end{proof}


\section{Subdivisions of \toric complexes}
\label{subdivisions}
In this section we introduce and study subdivisions of
\toric complexes. In parti\-cular, we will relate
the face rings of a \toric complex and its regular subdivisions.

\begin{defn}
Let $S$ and $T$ be a \toric complexes.
\begin{enumerate}
\item
$S$ is a {\em partial subdivision} of $T$ if
the following is satisfied:
\begin{enumerate}
\item  
$\Gen(S) \subseteq \Gen(T)$,
\item
for all $F\in \Pi_S$ there 
exists a face of $T$ containing $F$,
\item
$S_F(a)=T_{F'}(a)$ for all $F \in \Pi_S$ and $a \in F$ where $F' =
  \bigcap_{F \subseteq G, \ G \in \Pi_T} G$,
\item
for $F,G \in \Pi_S$ with $F \subseteq G$ 
we have $S_{F,G} = T_{F'G'}$.
\end{enumerate}
\item
$S$ is a {\em subdivision} of $T$ if $S$ is a partial subdivision of
$T$ satisfying that
$\cone(T_G(G))=\bigcup_{F\in \Pi_S, F\subseteq G} \cone(T_G(F))$ for
$G \in \Pi_T$. 
\item
$S$ is a {\em triangulation} of $T$ if $S$ is
simplicial and a subdivision of $T$.
\end{enumerate}
\end{defn}

In particular, if $S$ is a subdivision of $T$, then
$$
  ||T|| = \colim{G \in \Pi_T} \cone(T_G(G)) 
  = \colim{G \in \Pi_T} \bigcup_{F\in \Pi_S, F\subseteq G} \cone(T_G(F)) 
  \cong \colim{F \in
  \Pi_S} \cone(S_F(F)) = ||S||,
$$
and similarly we see that $|T| \cong |S|$.
\begin{ex}
  Given $F \subseteq \Z^d \setminus \{0\}$ and $a \in F$ we 
  construct an embedded toric complex $T(F,a)$ such that $T(F,a)$ is a
  subdivision of the toric complex $T(F)$ of Example \ref{first_ex}
  with $\Gen(T(F,a)) = \Gen(T(F))$. A subset $G$
  of $F$ is in $\Pi_{T(F,a)}$ if there exists a face $H$ of $F$ in
  $T(F)$ such that $a \notin H$ and $G = \cone(H \cup \{a\}) \cap
  F$ or $G=H$. We leave it as an instructive exercise for the reader
  to check 
  that $T(F,a)$ is a subdivison of $T(F)$.
\end{ex}
\begin{ex}
  Let $T$ be an embedded toric complex and let $a \in \Gen(T)$. For a
  face $F$ of $T$ with $a \in F$ we have defined the subdivision $T(F,a)$
  of $T(F)$. If $a \notin F$ we let $T(F,a) = T(F)$. The embedded
  toric complex 
  $S$ with $\Pi_S = \cup _{F \in \Pi_T} \Pi_{T(F,a)}$ is a subdivision
  of $T$. In the 
  case where $T=T(\Sigma,(G_{\sigma})_{\sigma \in \Sigma}, \varphi)$
  for a fan $\Sigma$ in $N$ as in Example \ref{fanex}, the subdivision
  $S$ of $T$ is well-known \cite[p. 48]{Fulton}.
\end{ex}

The following construction is inspired by the way
Sturmfels \cite{ST} subdivides the cone generated by a finite subset
$F$ of $\Z^d$. We have modified Sturmfels' construction slightly
so that the subdivision of a toric complex is again a 
toric complex.
\begin{defn}
\label{subdiv}
Let $T$ be a \toric complex and 
$\omega \colon \Gen(T)\to \R$.
Then the {\em $\omega$-subdivision} $\sd_\omega T$ of $T$ is 
given by the 
following data:
\begin{enumerate}
\item 
A subset $F$ of $\Gen(T)$ is
in $\Pi_{\sd_\omega T}$ if there 
exists a face of $T$ containing $F$ and letting 
$F' = \cap_{F \subseteq G \in \Pi_T} G$
there exists a linear form $\alpha_F$ on $\R^{d_{F'}}$ such that 
$F = \{a \in F' \colon \alpha_F({T_{F'}(a)}) = \omega(a)\}$ 
and 
$\alpha_F({T_{F'}(a)}) \leq {\omega(a)}$ for every $a \in F'$,
\item 
$(\sd_\omega T)_F(a) = T_{F'}(a)$ for $F \in
\Pi_{\sd_\omega T}$ 
and $F' = \cap_{F \subseteq G \in \Pi_T} G$,
\item
$(\sd_\omega T)_{FG} = T_{F'G'}$.
\end{enumerate}
If $\sd_\omega T$ is a subdivision (triangulation) of $T$, 
we say that $\sd_\omega T$ is a {\em regular subdivision (triangulation)} of $T$.
\end{defn}

A standard example of a non-regular triangulation is given in
\cite[Example 8.2]{ST}.
Let $e_1,\dots,e_{d+1}$ be the standard basis for $\R^{d+1}$.
A face $F$ of a cone $C \subset \R^{d+1}$ is called a {\em lower face}
if for all $x \in F$ and $r>0$ we have that $x-r\cdot e_{d+1} \not\in F$.
It is easy to see that $F$ is a lower face if and only if it
is the intersection of $C$ with a hyperplane $H_\alpha(0)$ such
that $\alpha(e_{d+1})<0$ and $C \subseteq H_\alpha^-(0)$.
\begin{rem}
\label{helper2}
Observe that (i) of Definition \ref{subdiv} 
can be reformulated as follows:
A subset $F$ of $\Gen(T)$ is in $\Pi_{\sd_\omega T}$ if 
$\cone(\sd_{\omega}T_{F'}(F))$ is the projection
of a lower face of 
$\cone( (T_{F'}(b),w(b)) \colon b\in F') \subset \R^{d_{F'} + 1}$
with respect to the last coordinate. 
\end{rem}

\begin{prop}
\label{sdomega}
$\sd_\omega T$ is a partial subdivision of $T$ for every
\toric complex $T$ and every
$\omega \colon \Gen(T) \to \R$. 
\end{prop}
\begin{proof}
(a), (b) and (c) of the definition of a \toric complex
are clearly satisfied.
It remains to verify conditions (i)--(iv).
(ii) and (iii) are fulfilled by definition.

(i):
Let $F, G \in \Pi_{\sd_{\omega} T}$.
Choose 
$F'=\cap_{F \subseteq F'' \in \Pi_T} F''$, 
$\alpha_F \colon \R^{d_{F'}} \to \R$
and
$G'=\cap_{G \subseteq G'' \in \Pi_{T}} G''$, 
$\alpha_G \colon  \R^{d_{G'}} \to \R$
as in Definition \ref{subdiv}.
$F' \cap G'$ is an element of $\Pi_T$. 
Via the inclusion of $\R^{ d_{F'\cap G'} }$ in $\R^{d_{F'}}$ 
induced by $T_{F'\cap G' F'}$ 
the linear form $\alpha_F$ on $\R^{d_{F'}}$
induces a linear form
$\beta_F$ on $\R^{d_{F'\cap G'}}$
such that
$\beta_F( T_{F'\cap G'}(a) )\leq \omega(a)$ for all $a \in F'\cap G'$
and
$$
F\cap G' = \{a \in F'\cap G'\colon \beta_F( T_{F'\cap G'}(a) )= \omega(a) \}.
$$
Analogously
there exists a
linear form
$\beta_G$ on  $\R^{d_{F'\cap G'}}$
such that
$\beta_G( T_{F'\cap G'}(a) )\leq \omega(a)$ for all $a \in F'\cap G'$
and
$$
F'\cap G = \{a \in F'\cap G'\colon \beta_G( T_{F'\cap G'}(a) )= \omega(a) \}.
$$
Let
$\alpha_{F \cap G} := (\beta_F + \beta_G)/2 \colon \R^{d_{F'\cap G'}}
\to \R.$
Then
$$
\alpha_{F \cap G}( T_{F'\cap G'}(a) )\leq  \omega(a) \text{ for all } a \in F'\cap G'
$$
and
$$
F\cap G = \{a \in F'\cap G'\colon \alpha_{F \cap G}( T_{F'\cap G'}(a)
)= \omega(a) \}. 
$$
Hence $F\cap G \in \Pi_{\sd_\omega T}$ and this shows (i).

(iv): 
Assume that $F \subseteq G \in \Pi_{\sd_\omega T}$.
If $F \in \Pi_{\sd_\omega T}$, then choose 
$F'$, $\alpha_F \colon \R^{d_{F'}} \to \R$ and
$G'$, $\alpha_G \colon \R^{d_{G'}} \to \R$ 
as in the proof of (i).
Note that $F' \subseteq G'$.
Therefore 
$\cone(T_{G'}(F'))=H_\gamma(0)\cap \cone(T_{G'}(G'))$ is a face of
$\cone(T_{G'}(G'))$ for some linear form $\gamma$ on $\R^{d_{G'}}$,
$\cone(T_{G'}(G')) \subseteq H^-_\gamma(0)$
and
$T_{G'}(F')=\cone(T_{G'}(F'))\cap T_{G'}(G')$.
We have the inclusion 
$\R^{d_{F'}}\to \R^{d_{G'}}$ 
and we extend $\alpha_F$ arbitrarily  to a linear form on $\R^{d_{G'}}$.

Now choose $t\gg 0$ such that for
$\beta_F = \alpha_F + t\cdot \gamma$ we have that
$F = \{a \in G' \colon \beta_F({T_{G'}(a)}) = \omega(a)\}$ 
and 
$\beta_F({T_{G'}(a)}) \leq {\omega(a)}$ for every $a \in G'$.
Thus $\beta=\beta_F- \alpha_G$ is a linear form
such that 
$\cone(T_{G'}(F))=H_\beta(0)\cap \cone(T_{G'}(G))$ 
is a face of $\cone(T_{G'}(G))$
and
$T_{G'}(F)=\cone(T_{G'}(F))\cap T_{G'}(G)$.

Conversely, assume that 
$F\subseteq G \in \Pi_{\sd_\omega T}$,
$G'$, $\alpha_G \colon \R^{d_{G'}} \to \R$ are chosen as above
and 
$\cone(T_{G'}(F))=H_\beta(0)\cap \cone(T_{G'}(G))$ 
is a face of $\cone(T_{G'}(G))$
defined by some linear form $\beta$ on $\R^{d_{G'}}$
and
$T_{G'}(F)=\cone(T_{G'}(F))\cap T_{G'}(G)$.
Define 
$\beta_F=\alpha_G+t\cdot\beta$. 
For a suitable $t\geq 0$ we have that
$F = \{a \in G' \colon \beta_F({T_{G'}(a)}) = \omega(a)\}$ 
and 
$\beta_F({T_{G'}(a)}) \leq {\omega(a)}$ for every $a \in G'$.
This implies that $F \in \Pi_{\sd_\omega T}$ and we are
also done in this case.

Hence $\sd_\omega T$ is a \toric complex
which is a partial subdivision of $T$ by construction.
\end{proof}

If $S$ is a subcomplex of a \toric complex $T$,
then a map 
$\omega \colon \Gen(T) \to \R$
induces a map
$\omega \colon \Gen(S) \to \R$.
\begin{cor}
\label{sd_cor}
Let $T$ be a \toric complex, 
$\omega \colon \Gen(T) \to \R$ 
and $S$ a subcomplex of $T$.
Then $\sd_\omega S$ is a subcomplex
of $\sd_\omega T$ and if
$\sd_\omega T$ is a regular subdivision (triangulation) of $T$, 
then $\sd_\omega S$ is a regular subdivision (triangulation) of $S$.
\end{cor}
\begin{proof}
This follows from Definition \ref{subdiv}
since the property to be a regular subdivision is 
defined on the faces of a \toric complex.
\end{proof}

A subcomplex of a \toric complex inherits regular subdivisions 
as noted in \ref{sd_cor}.
The corresponding face rings are related as follows. 
\begin{cor}
\label{main_kor2}
Let $T$ be a \toric complex, $\omega \colon \Gen(T) \to \R$ 
such that
$\sd_\omega T$ is a regular subdivision of $T$
and $S$ a subcomplex of $T$.
Then
$$
K[\sd_\omega S] 
\cong 
K[\sd_\omega T]/
(x_F \colon F\subseteq \Gen(\sd_\omega T) \text{ is not contained in any } G \in \Pi_{\sd_\omega S}).
$$
\end{cor}
\begin{proof}
By Corollary \ref{sd_cor} 
the \toric complex 
$\sd_\omega S$ is a regular subdivision of
$S$ and a subcomplex of $\sd_\omega T$.
The isomorphisms follow from Proposition \ref{toric_proj}.
\end{proof}

The toric complex $\sd_\omega T$ is not always a subdivision of $T$. 
As an illustration let $F$ denote the finite set $F = \{-1,1\}
\subseteq \Z$, let $\omega \colon F \to \R$ be the constant function with
value $-1$ and consider
the toric complex $T(F)$ of Example \ref{first_ex}. Then $||T(F)|| = \R$ and
$||\sd_\omega T(F)|| = \emptyset$, 
because $\Pi_{\sd_\omega T(F)}$ is the empty set.
\begin{prop}
\label{criteria}
Let $T$ be a \toric complex and 
$\omega \colon \Gen(T) \to \R$.
Then $\sd_\omega T$ is a regular subdivision of $T$ 
if one of the following conditions is satisfied:
\begin{enumerate}
  \item
  All values $\omega(a)$ of $\omega$ are positive.
  \item
  $T$ is an $\N$-graded \toric complex.
\end{enumerate}
\end{prop}
\begin{proof}
In both cases it remains to show that
$$
\cone(T_G(G))=\bigcup_{F\in \Pi_{\sd_\omega T}, F\subseteq G} \cone(T_G(F))
$$ 
for
$G \in \Pi_T$, because by \ref{sdomega}
we know already that $\sd_\omega T$ is a partial subdivision of $T$.
Let $G \in \Pi_{T}$ and $0 \ne x\in\cone(T_{G}(G))$. 
Consider the cone 
$$
C=\cone((T_G(b),\omega(b))\colon b\in G) \subset \R^{d_G +1}
\text{ and }
P=\{(x,t)\colon t \in \R  \} \subset \R^{d_G +1}.
$$
Since $x \in \cone(T_{G}(G))$
we have that $C\cap P \neq \emptyset$.
If one of the conditions (i) or (ii) is satisfied
there exists
$$
s =\inf \{t \in \R \colon (x,t) \in P\cap C\} >-\infty.
$$
Then $(x,s)$ is an element of a lower face $F$ of $C$
and by Remark \ref{helper2} 
there exists
a linear form $\alpha$ on $\R^{d_G}$
and $F'\in \Pi_T$ such that
$\alpha_F({T_{F'}(a)}) \leq {\omega(a)}$ for every $a \in F'$
with $x\in \cone(F)$ and
$F = \{a \in F' \colon \alpha_F({T_{F'}(a)}) = {\omega(a)})\}$.
\end{proof}
\begin{ex}
  Let $T$ be a toric complex and consider the function $\omega \colon
  \Gen(2T) \to \R$ with $\omega(2\chi_a) = 2$ for $a \in \Gen(T)$ and
  $\omega(u) = 1$ if $u$ is not of the form $2\chi_a$. Here $\chi_a$ is
  the indicator function on $a$ with $\chi_a(a) = 1$ and
  $\chi_a(b) = 0$ for $b \in \Gen(T) \setminus \{a\}$. 
  In the
  case, where $T = T(\Sigma,(G_{\sigma})_{\sigma \in \Sigma},
  \varphi)$ for a rational fan $\Sigma$, the toric complex $2T$ is of
  the form $2T = T(\Sigma,(H_{\sigma})_{\sigma \in \Sigma},
  \varphi)$, where $H_\sigma = \{a + b \colon a,b \in
  G_{\sigma}\}$ and the fan $\{\cone(T_F(F)) \colon F \in
  \Pi_{sd_{\omega}(2T)} \}$ is a subdivision of the fan $\Sigma$.
\end{ex}

Next we construct a homomorphism
$$
\varphi\colon
K[\Gen(T)]/
\init_\omega(I_T) \to K[\sd_\omega T] \cong 
\prolim{F \in \Pi^{\op}_{\sd_\omega T}} K[T(F)].
$$ 
If we show that for all $F \in \Pi_{\sd_\omega T}$ the projection
$
i^* \colon K[\Gen(T)] \to K[F]$ induced by the inclusion $F \subseteq \Gen(T)$
induces a homomorphism 
$
\varphi_F \colon K[\Gen(T)]/\init_\omega(I_T) \to K[T(F)],
$
then $\varphi$ exists by the universal property of the limit.

Given $F \in \Pi_{\sd_\omega T}$ we let 
$F' = \cap_{F \subseteq G \in \Pi_T} G \in \Pi_T$ as above. 
By the construction of $\sd_\omega T$, 
there exists a linear form $\alpha_{F}$ on $\R^{d_{F'}}$
such that 
$F = \{ a \in F' \colon \alpha_F(T_{F'}(a))  = \omega (a)\}$ 
and 
$\alpha_{F}(T_{F'}(a)) \leq \omega(a)$ for $a \in F'$. 
Note that $\omega$ induces a weight order on $K[F']$. 
For $x^u \in  K[F']$ 
we have that 
$$
\omega(u)=
\sum_{a \in F'} u(a)\omega(a) \ge \sum_{a \in F'} u(a) \alpha_F(T_{F'}(a))
=
\alpha_F(T_{F'}(u))
,
$$ 
and equality
holds precisely if $\supp(u) \subseteq F$.
Let $x^v \in K[F']$ with
$
T_{F'}(u) 
= 
T_{F'}(v).
$
If 
$$
\sum_{a \in F'} u(a) \omega(a) = \sum_{a \in F'} v(a) \omega(a),
$$ 
then $\supp(u) \subseteq F$ if and only if $\supp(v) \subseteq F$. 

Let 
$g \in \Ic_{F'} \subset K[F']$.
We show that the projection of
$\ini_{\omega}(g)$ to $K[F]$ is an element of $\Ic_{F}$.  
The ring $K[F']$ is $\Z^{d_{F'}}$-graded by giving $x_a$ the degree $T_{F'}(a)$ for $a\in F'$ and
$\Ic_{F'}$ is homogeneous with respect to this grading. 
Thus 
without loss of generality we may assume that
there exists $z\in \Z^{d_{F'}}$ such that
$g=\sum_{u}c_ux^u$ and $ T_{F'}(u)=z$ for  $c_u\neq 0$. 

By the above discussion
either $\ini_{\omega}(g)=g$ 
or for all monomials $x^u$ in $\ini_{\omega}(g)$ we have
$\supp(u)\not\subseteq F$. 
In the first case the image of $\ini_{\omega}(g)=g$ 
under the projection $i^* \colon K[F'] \to K[F]$ is in $\Ic_F$.
In the second case $i^*(\ini_{\omega}(g)) = 0$.
We obtain that  
$i^* K[F'] \to
K[F]$ induces a homomorphism 
$\varphi_{FF'} \colon K[F']/\init_\omega (\Ic_{F'}) \to K[F]/\Ic_F.$ 
The following diagram of natural projections 
commutes for $F \in \Pi_{\sd_\omega T}$:
\begin{displaymath}
  \begin{CD}
  K[\Gen(T)] @>>> K[\Gen(\sd_\omega T)] @>>> 
  K[\Gen(\sd_\omega T)]/I_{\sd_\omega T} \\
  @VVV @. @VVV \\
  K[\Gen(T)]/\init_{\omega}(I_T) @>>> 
  K[F']/\init_\omega(\Ic_{F'}) @>{\varphi_{FF'}}>> K[F]/\Ic_F.
  \end{CD}
\end{displaymath}
Taking limit with respect to $F$  we get a
commutative diagram of the form: 
\begin{displaymath}
  \begin{CD}
  K[\Gen(T)] @>>> K[\Gen(\sd_\omega T)] @>>> 
  K[\Gen(\sd_\omega T)]/I_{\sd_\omega T} \\
  @VVV @. @V{\cong}VV \\
  K[\Gen(T)]/\init_{\omega}(I_T) @= 
  K[\Gen(T)]/\init_{\omega}(I_T)
  @>>> \prolim{F \in \Pi^{\op}_{\sd_\omega T}} K[F]/\Ic_F.
  \end{CD}
\end{displaymath}
Thus the projection 
$K[\Gen(T)] \to K[\Gen(\sd_\omega T)] $
induces a homomorphism
$\varphi$
from $K[\Gen(T)]/
\init_\omega(I_T)$ to $K[\sd_\omega T]$. 
Since 
$K[\sd_\omega T]$ is reduced 
we obtain a homomorphism 
$
\Phi
\colon
K[\Gen(T)]/\rad(\init_{\omega}(I_T)) \to K[\sd_\omega T].
$

\begin{thm}
\label{main_result}
For every \toric complex $T$ 
and 
$\omega \colon \Gen(T) \to \R$ 
such that  
$\sd_\omega T$ is a regular subdivision of $T$ the homomorphism
$
\Phi \colon K[\Gen(T)]/
\rad(\init_\omega (I_T)) \to
K[\sd_\omega T] 
$ 
is an isomorphism. 
\end{thm}
\begin{proof}
Recall that 
$I_{\sd_\omega T} =  J_{\sd_\omega T} + \sum_{F \in \Pi_{\sd_\omega T}} I_F$
and 
$I_T =  J_T + \sum_{G \in \Pi_T} I_G$. 

We start by showing that the inclusion 
$i_* \colon K[\Gen(\sd_\omega T)] \to K[\Gen(T)]$
maps $I_{\sd_\omega T}$ into 
$\init_\omega(I_T)$. 
We identify monomials in
$K[\Gen(\sd_\omega T)]$ with the corresponding ones in $K[\Gen(T)]$.
For every $F \in \Pi_{\sd_\omega T}$ there
exists $F' \in \Pi_T$ with $F \subseteq F'$, and 
a linear from $\alpha_F$ on $\R^{d_{F'}}$ 
with $\alpha_F(T_{F'}(a)) =\omega(a)$ for
$a \in F$ and 
$\alpha_F(T_{F'}(a)) < \omega(a)$ for $a \in F' \setminus F$. 
The (image of a) binomial 
$x^u-x^v \in \Ic_F$ is an element of $I_{F'}$. 
(Recall that these binomials generate $I_F$.)
Since $\supp(u),\supp(v) \subseteq F$
we have
$$
\omega(u) 
= \sum_{a \in F} u(a) \omega(a) 
= \sum_{a \in F} u(a) \alpha_F(T_{F'}(a)) 
= \alpha_F(T_{F'}(u)) 
= \alpha_F(T_{F'}(v)) 
$$
$$
= \sum_{a \in F} v(a) \alpha_F(T_{F'}(a)) 
= \sum_{a \in F} v(a) \omega(a) 
= \omega(v),
$$ 
thus $x^u - x^v$ is an element of
$\init_{\omega}(I_{F'}) \subseteq \rad(\init_\omega (I_T)$). 

If $x^u \in J_{\sd_\omega T}$, 
then either $x^u$ is in $J_T$ and we are done, 
or there exists an 
$G \in \Pi_T$ such that $\supp(x^u) \subseteq G$. 
Suppose that $x^u \in J_{\sd_\omega T}$ and that  
$\supp(x^u) \subseteq G \in \Pi_T$.  
Since $\sd_\omega T$ is a subdivision of $T$ and 
$T_{G}(u)  \in \cone(T_{G}(G))$ 
there exists $F \in \Pi_{\sd_\omega T}$
with $F \subseteq G$ such that $T_G(u) \in \cone(T_G(F))$. 
We may assume that $G=F'$ and then
there exists a linear from
$\alpha_F$ on $\R^{d_{F'}}$ 
with $\alpha_F(T_{F'}(a)) =\omega(a)$ for
$a \in F$ and 
$\alpha_F(T_{F'}(a)) < \omega(a)$ for $a \in F' \setminus F$.

Since $T_{F'}(u) \in \cone(T_{F'}(F))$ 
there exists a map $\lambda\colon F \to \R_+$ 
such that $T_{F'}(u) = \sum_{a \in F} \lambda(a) T_{F'}(a)$. 
Using for example that Farkas lemma holds over both $\R$ and $\Q$, 
we may assume that $\lambda$ takes values in $\Q_+$, 
and thus there exists 
$v \colon F \to \N$ and $n \in \N$
such that $T_{F'}(nu) = T_{F'}(v)$. 
Now consider
$x^{nu} - x^v \in I_{F'}$. 
There
exists at least one 
$a \in \supp(u)$ with $\alpha_F(T_{F'}(a)) < \omega(a)$
since $x^u \in J_{\sd_\omega T}$. 
Since $\supp(v) \subseteq F$ we have
$$
\omega({nu}) > \alpha_F(T_{F'}(nu)) 
= \alpha_F(T_{F'}(v)) 
= \omega(v).
$$ 
We conclude that 
$x^{nu} \in \init_\omega (I_{F'}) \subseteq  \init_\omega (I_T)$.
Hence
$x^u \in \rad(\init_\omega (I_T))$. 
This finishes the proof of the fact that the inclusion 
$i_* \colon K[\Gen(\sd_\omega T)] \to K[\Gen(T)]$ 
maps $I_{\sd_\omega T}$ into $\rad(\init_\omega (I_T))$. 

We denote 
by 
$
\Psi \colon K[\Gen(\sd_\omega T)]/I_{\sd_\omega T} \to 
K[\Gen(T)]/\rad(\init_\omega (I_T))
$
the map induced by  
$i_*$.
It is immediate that $\Phi \circ\Psi=\id_{K[\sd_\omega T]}$. 
We need to prove
that $\Psi\circ \Phi$ is the identity, or equivalently that $\Psi$ is
onto. Since $\sd_\omega T$ is a subdivision of $T$ there exists for every $G
\in \Pi_T$ and $a \in G$ an 
$F \subseteq G$ with $F \in \Pi_{\sd_\omega T}$ such that 
$T_G(a) \in \cone(T_G(F))$.
Again we may assume that $G=F'$ and $\alpha_{F'} \colon \R^{d_{F'}}
\to \R$
is chosen as above.
If $x_a$ is not in the image of the inclusion 
$i_* \colon K[\Gen(\sd_\omega T)] \to K[\Gen(T)]$, 
then $a$ is not in $F$ and 
$\alpha_F(T_{F'}(a)) < \omega(a)$. 
Writing $n\cdot T_{F'}(a) = \sum_{b \in F} v(b) T_{F'}(b)$ 
for $n \in \N$ and $v \colon F \to \N$ as above, 
we compute that
$x_a^{n} - x^v \in I_{F'}$, and 
$\omega(na) > \alpha_{F}(n\cdot T_{F'}(a)) =\alpha_F(T_{F'}(v)) 
= \omega(v)$. 
Hence $x_a^{n} \in \init_{\omega} (I_{F'})$. 
It follows that 
$x_a \in \rad(\init_{\omega} (I_{F'})) \subseteq \rad(\init_{\omega} (I_{T}))$ 
for every $a \in \Gen(T)  \setminus \Gen(\sd_\omega T)$. 
This finishes the proof of the theorem.
\end{proof}
Theorem \ref{main_result} generalizes a result of Sturmfels
\cite[Theorem 8.3]{ST}: 
\begin{cor}
\label{main_result2}
Let $T$ be a \toric complex and $\omega \colon \Gen(T) \to \R$
such that $\sd_\omega T$ is a regular subdivision of $T$. Then
\begin{enumerate}
  \item
  $\rad(\ini_{\omega} (I_T))$ is a square-free monomial ideal if and
  only if 
  $\sd_{\omega} T$ is a regular triangulation of $T$.
  In this case the abstract simplicial complex induced by
  $\rad(\ini_{\omega} (I_T))$ coincides with $\Pi_{\sd_{\omega} T}$.
  \item
  Given $\omega'
  \colon \Gen(T) \to \R$ such that $\sd_{\omega'} T$ is a subdivision of
  $T$ we have that $\sd_\omega T = \sd_{\omega'} T$ if and only if 
  $\rad(\ini_{\omega} (I_T)) = \rad(\ini_{\omega'} (I_T))$.
\end{enumerate}
\end{cor}
\begin{proof}
(i):
By Theorem \ref{main_result} $I_{\sd_\omega T}$ is a monomial ideal
if and only if $\rad(\ini_{\omega} (I_T))$ is a square-free monomial ideal.
Hence it follows from \ref{helper3} that
$\rad(\ini_{\omega} (I_T))$ is a square-free monomial ideal if and only if
$\sd_{\omega} T$ is a regular triangulation of $T$ and
that
the abstract simplicial complex induced by
$\rad(\ini_{\omega} (I_T))$ coincides with $\Pi_{\sd_{\omega} T}$.

(ii):
We only show that 
$\rad(\ini_{\omega} (I_T)) = \rad(\ini_{\omega'} (I_T))$ implies that
$\sd_\omega T = \sd_{\omega'} T$ since the other implication is a
direct consequence of \ref{main_result}.
Observe that by \ref{main_result} we have that 
$\rad(\ini_{\omega} (I_T)) = \rad(\ini_{\omega'} (I_T))$ if and only
if $\Gen(\sd_\omega T) = \Gen(\sd_{\omega '}T)$ and 
$I_{\sd_\omega T} = I_{\sd_{\omega '} T}$.
The result follows from the fact that for a monomial $x^u$ the
support of $u$ is contained in a face of $\sd_\omega T$ if and only if
$x^u \notin I_{\sd_{\omega} T}$.
\end{proof}

Recall that given a subset $L$ of $\Z^d$ 
we let $M_L$ denote the submonoid of $\Z^d$ generated by $L$.
Let $L \subseteq L'$ be subsets of $\Z^{d}$. 
We say that $M_L$ is {\em integrally closed} in $M_{L'}$ 
if for all $x \in M_{L'}$ with $n\cdot x \in M_L$
for some $n \in \N$ we have that $x \in M_L$.
It is easy to see that this is equivalent to the fact that
$\cone(L) \cap M_{L'} = M_{L}$.

The following proposition can be 
extracted from the proof of 
Theorem \ref{main_result}: 
\begin{prop}
\label{radi_prop}
Let $T$ be a toric complex and $\omega\colon \Gen(T) \to \R$. 
If $\Gen(\sd_{\omega} T) = \Gen(T)$ and 
$T_G(F)$ is integrally closed in $T_G(G)$ 
for all 
$F \in \Pi_{\sd_\omega T}$ and $G \in \Pi_T$ with $F \subseteq G$, 
then the ideal 
$\ini_\omega(I_T)$ is a radical ideal.
\end{prop}


\section{Edgewise subdivision of a simplicial complex} 
\label{edgewisesubdivision}
In this section $\Delta$ denotes an abstract 
simplicial complex on the vertex set $V = \{1,\dots,d+1\}$,
i.e.\ $\Delta$ is a set of subsets of $V$ 
and $F\subseteq G$ for $G \in \Delta$ implies that $F \in \Delta$.
Let $K[\Delta]=K[V]/I_\Delta$ be the Stanley-Reisner ring of $\Delta$ 
where 
$I_\Delta=(x_F \colon F \subset \{1,\dots,d+1\}, F \not\in\Delta)$
is the Stanley-Reisner ideal of $\Delta$.
Observe that $K[\Delta]=K[T(\Delta)]$
where $T(\Delta)$ is the \toric complex
of \ref{second_ex}.
By Theorem \ref{toric_veronese} the ring
$K[rT(\Delta)] =
K[T(\Delta)]^{(r)}$ is the $r$-th Veronese subring of $K[T(\Delta)]$.
Given an element $m$ of $\Z^{d+1}$ we let $m_j$
denote its $j$-th coordinate for $1 \le j \le d+1$.

Guided by the regular subdivision introduced by 
Knudsen and Mumford \cite[pp. 117--123]{Mumford} 
and
the edgewise subdivision of simplicial sets  
(see for example \cite{BHM}, \cite{EDGR} or \cite{Gr})
we define the following simplicial complex.

\begin{defn}
\label{edgewise_subdiv}
Let $\Delta$ be an abstract simplicial complex with vertex set $V$, let
$r \ge 1$ and $\iota_V \colon V \to \Z^{d+1}$ as in 
Example \ref{second_ex}. 
Let $\iota_V$ also denote the map 
$\N^V \to \Z^{d+1}, u \mapsto \sum_{i=1}^{d+1} u(i)\iota_V(i)$.
The {\em 
$r$-fold edgewise subdivision} of 
$\Delta$ is the abstract simplicial complex 
$\esd_r (\Delta)$ consisting of the subsets
$F \subseteq rV = \{ u \in \N^V \colon \sum_{i \in V} u(i) = r \}$ with
$\bigcup_{u \in F} \supp(u) \in \Delta$ 
satisfying that for $u,u' \in F$ we either have that 
$0 \le (\iota_V(u) - \iota_V(u'))_j \le 1$ for every
$j \in \{1,\dots,d+1\}$ or that
$0 \le (\iota_V(u') - \iota_V(u))_j \le 1$ for every
$j \in \{1,\dots,d+1\}$.
\end{defn}
The name
$r$-fold edgewise subdivision comes from the fact that
the edges of $\Delta$ are subdivided in $r$ pieces.
The goal of this section will be to
relate $K[\esd_r(\Delta)]$ to the $r$-th
Veronese subring of the Stanley-Reisner ring $K[\Delta]$.

Observe that 
we
have required that 
the coordinate-wise partial order on $\Z^{d+1}$ induces a total order
on $\iota_V(F)$ 
for every $F \in \esd_r(\Delta)$. If $m(F)$ denotes the minimal element of
$\iota_V(F)$ in this order, then we also require that
the set $\iota_V(F) \subseteq \Z^{d+1}$
is a subset of $m(F) + (\{0,1\}^{d}\times \{0\})$, the vertices of a  
$d$-dimensional lattice cube.

Given an element $\sigma$ of the symmetric group $\Sigma_d$ 
let 
\begin{displaymath}
  \Delta_{\sigma} = \conv (0,e_{\sigma(1)},e_{\sigma(1)} +
  e_{\sigma(2)}, \dots, e_{\sigma(1)} +
  e_{\sigma(2)} +\dots + e_{\sigma(d)} ) \subseteq \R^{d+1},
\end{displaymath}
and denote by $\Delta^{d}$ the abstract simplex on the vertex set
$V$.
For every maximal face $F$ of $\esd_r(\Delta^d)$
the polytope $\conv(\iota_V(F))$ is a simplex of the form $m(F) +
\Delta_{\sigma(F)}$ for a unique $\sigma(F) \in \Sigma_d$. Conversely, for
every $m \in \Z^{d+1}$ and $\sigma \in \Sigma_d$ such that $m +
\Delta_\sigma$ is contained in
$\conv(\iota_V(rV)) = r \conv(\iota_V(V))$, 
the set $F$ corresponding to the vertices of $m +
\Delta_\sigma$ is a maximal element of $\esd_r(\Delta^d)$. 

Given $k$ with $0 < k < r$ and $y \in \R^{d+1}$
we define:
\begin{displaymath}
  \alpha^{ij}_k(y) = 
  \begin{cases}
    r(y_i - y_j) + k y_{d+1} & \text{for $1 \le i < j < d+1$}, \\
    ry_{i} - ky_{d+1}        & \text{for $1 \le i < j = d+1$}.
  \end{cases}
\end{displaymath}
The following easy and well-known result implies that
the linear form $\alpha^{ij}_k \colon \R^{d+1} \to \R$
has the property that $\conv(\iota_V(F)) =
m(F) +
\Delta_{\sigma(F)}$ is contained in either its positive- or its negative
associated 
halfspace for every $F \in \esd_r(\Delta^d)$. 
For a proof we refer to Knudsen \cite[Lemma 2.4]{Mumford}.
\begin{lem}
\label{reg_sdiv}
  The
  hyperplanes $H^{ij}_k =
  (\alpha^{ij}_k)^{-1}(0)$ induce a triangulation of the simplex
  $\conv(\iota_V(rV))$ into 
  simplices of the form $m + \Delta_\sigma$. 
\end{lem}
The above lemma in particular implies that after choosing a total
order on $rV$ we have that 
$|T(\esd_r(\Delta^d))| \cong
\conv(\iota_V(rV)) \cong \Delta_{e}$ where $e \in \Sigma_d$ is the neutral
element.    

\begin{defn}
\label{esd_toric}
  Let $\Delta$ be an abstract simplicial complex and let $r \ge
  1$. Define $\esd_r (T(\Delta))$ to be the
  (embedded) \toric complex with $\Pi_{\esd_r (T(\Delta))} = \esd_r(\Delta)$,
  with
  $\esd_r(T(\Delta))_F$ given by the restriction of
  $\iota_V \colon rV \to \Z^{d+1}$ to $F$
  for $F \in \Pi_{\esd_r(T(\Delta))}$ 
  and with
  $\esd_r(T(\Delta))_{FG} = \id_{\Z^{d+1}}$ for $F,G \in
  \Pi_{\esd_r(T(\Delta))}$ with $F \subseteq G$.
\end{defn}
Define the convex function $f \colon \R^{d+1} \to \R$ by 
$$
f(y) = \sum_{\stackrel{1 \le i < j
\le d+1} {0<k<r}}  |\alpha^{ij}_k(\iota_V(y))|
$$
and let $\omega = f \circ \iota_V \colon rV \to \R$.
Lemma \ref{reg_sdiv} implies that for every $F \in
\esd_r(\Delta^d)$ there exists a linear form 
$\alpha_F \colon \R^{d+1} \to \R$ such
that 
$\alpha_F(y) = f(y)$ for every $y \in m(F) + \Delta_{\sigma(F)}$.
Since the convex function $f$ agrees with $\alpha_F$ on an open
subset of $\R^{d+1}$ we have that 
$\alpha_F(y) \le f(y)$ for $y \in \R^{d+1}$, and because $m(F)
+ \Delta_{\sigma(F)}$ is a simplex in the triangulation induced by the 
$H^{ij}_k$, equality
holds precisely if $y \in \conv(\esd_r(T(\Delta^d))_F(F))$. This implies that 
$\esd_r(\Delta^d) \subseteq \Pi_{\sd_\omega rT(\Delta^d)}$.
On the other hand if 
$F$ is a maximal face of ${\sd_\omega(rT(\Delta^d))}$, then 
by \ref{reg_sdiv}
$\conv(\esd_r(T(\Delta^d))_F(F))$ is of the form $m(F) + 
\Delta_{\sigma(F)}$. 
As in the discussion after \ref{edgewise_subdiv} we see 
that $F \in
\esd_r(\Delta^d)$. 
\begin{prop}
  For every abstract simplicial complex $\Delta$ with vertex set $V$ we have
  that $\sd_\omega rT(\Delta) = \esd_r(T(\Delta))$ is a regular
  triangulation of $rT(\Delta)$.
\end{prop}
\begin{proof}
  We have just proved the result in the case 
  $\Delta =\Delta^d$. Since 
  $\Pi_{\sd_\omega rT(\Delta)} = 
  \Pi_{\sd_\omega rT(\Delta^d)} \cap \{ F \subseteq rV \colon \bigcup_{u \in F} \supp(u) \in \Delta\}$ 
  and 
  $\Pi_{\esd_r T(\Delta)} = \Pi_{\esd_r(T(\Delta^d))}\cap  
  \{ F \subseteq rV \colon \bigcup_{u \in F} \supp(u) \in \Delta\}$ 
  it follows
  that $\sd_\omega rT(\Delta)
  = \esd_r(T(\Delta))$. Applying Corollary \ref{sd_cor} we get that $\sd_\omega
  rT(\Delta)$ is a 
  regular triangulation of $rT(\Delta)$. 
\end{proof}
\begin{cor}
\label{inesd}
We have that
$$
K[\esd_r(\Delta)] =
K[\esd_r(T(\Delta))] = K[\sd_\omega(rT(\Delta))] \cong
K[rV]/\ini_\omega(I_{rT(\Delta)}).
$$
\end{cor}
\begin{proof}
  Only the last isomorphism does not follow directly from the
  definitions. By Theorem \ref{main_result} 
  the face ring $K[\sd_\omega(rT(\Delta))]$ is isomorphic to
  $K[rV]/\rad(\ini_\omega(I_{rT(\Delta)}))$.
  Since for every $F \in \Pi_{\sd_\omega rT(\Delta)}$ 
  the set $\iota_V(F)$ is a subset of the
  vertices of a simplex of the form $m + \Delta_\sigma$ for $m \in
  \Z^d \times \{r\}$ and $\sigma \in \Sigma_d$, the set $\iota_V(F)$
  can be extended to a basis for $\Z^{d+1}$. By Proposition
  \ref{radi_prop} we conclude that 
  $\ini_\omega(I_{rT(\Delta)})$ is a radical ideal.
\end{proof}
Let $\psi \colon V^r \to rV$ be the surjective map
taking $v=(v_1,\dots,v_r)$ to the function
$\psi(v)\colon V \to \N$ 
with $\psi(v)(i)=|\{l \in \{1,\dots,r\}: v_l=i  \}|$.
Restricting $\psi$ to the subset $W$ of
$V^r$ consisting of tuples of the form $v = (v_1,\dots,v_r)$ with 
$1\le v_1 \le v_2 \le \dots \le v_{r} \le d+1$ we obtain a bijection
$\psi \colon W \to rV$. 
Note that $(\iota_V(\psi(v)))_j$ is given by the cardinality of
$\{l\in\{1,\dots,r\}:1 \leq v_l\leq j \}$ for $1 \le j \le d+1$. 
Given a subset $F$ of $W$ we have that
$\psi(F) \in \esd_r(\Delta)$ if and only if 
the following two
conditions are satisfied:
\begin{enumerate}
\item
there exists an ordering
$F = \{v_1, \dots, v_s\}$ of the elements of $F$ with  
$v_i =(v_{i1}, \dots,v_{ir}) \in W$
such that
$$
\qquad \quad \, 1 \le v_{11}\leq v_{21} \leq \cdots \leq v_{s1} \leq 
v_{12}\leq v_{22} \leq \cdots  \leq
v_{1r}\leq v_{2r} \leq \cdots \leq v_{sr} \le d+1,
$$
\item
the set $\{v_{ij} \colon 1\leq i \leq s,\ 1\leq j \leq r  \}$
is a face of $\Delta$.
\end{enumerate}

Let $K[W]$ denote the polynomial ring on the set $W$.
We have a natural surjection
$
\Phi_{rT(\Delta)}
\colon K[W] \to 
K[\Delta]^{(r)}
$
defined by
$\Phi_{rT(\Delta)} 
(x_{(v_1,\ldots,v_r)}) = x_{v_1}\cdots x_{v_r}
$
and we write 
$I'_{rT(\Delta)}$ for the kernel of $\Phi_{rT(\Delta)}$.
Note that $\Phi_{rT(\Delta)}$ induces an isomorphism
$
\Phi_{rT(\Delta)} \colon K[W]/I'_{rT(\Delta)} \to
K[\Delta]^{(r)}=K[rT(\Delta)]=K[rV]/I_{rT(\Delta)}$.   

A monomial $x^u \in K[W]$ is of the form
$x^u=\prod_{i=1}^sx_{u_i}$ for $u_i=(u_{i1},\dots,u_{ir})\in W$. 
There exists a unique 
matrix of the form 
\begin{equation*}
\label{star}
\begin{pmatrix}
v_{11} & v_{12} & \cdots & v_{1r}\\
v_{21} & v_{22} & \cdots & v_{2r}\\
\vdots & \vdots & \vdots & \vdots\\
v_{s1} & v_{s2} & \cdots & v_{sr}
\end{pmatrix}
\end{equation*}
with
$$
1 \le v_{11}\leq v_{21} \leq \cdots \leq v_{s1} \leq 
v_{12}\leq v_{22} \leq \cdots  \leq
v_{1r}\leq v_{2r} \leq \cdots \leq v_{sr} \le d+1.
$$
and 
$\prod_{i=1}^s\prod_{j=1}^r x_{u_{ij}}=
\prod_{i=1}^s\prod_{j=1}^r x_{v_{ij}} \in K[V]$.
Let $\sort(x^u) = \prod_{i=1}^s x_{(v_{i1},\dots,v_{ir})}$. 
Motivated by Sturmfels \cite[14.2]{ST}
we call the
monomial $x^u \in K[W]$ {\it sorted} if $x^u = \sort(x^u)$. 

If $x^u  = x_{(u_{11},\ldots,u_{1r})}x_{(u_{21},\ldots,u_{2r})}$, 
then $\psi(u_{11},\ldots,u_{1r})$ and $\psi(u_{21},\ldots,u_{2r})$
are connected by an edge in $\esd_r(T(\Delta^d))$ if and only if $x^u$
is sorted.

Let $\omega' = \omega \circ \psi \colon W \to \R$
where $\omega=f \circ \iota_V$ as above. 
The following is essentially
Theorem 14.2 in Sturmfels \cite{ST}
(see Hibi-Ohsugi \cite{OHHI} for related Gr\"obner bases).

\begin{prop}
\label{first_step}
The initial ideal $\ini_{\omega'} (I'_{rT(\Delta^{d})})$ 
is generated by the initial polynomials of 
the set $\Gc$ consisting of the binomials
$$
x_{(u_{11},\ldots,u_{1r})}x_{(u_{21},\ldots,u_{2r})}
-
\sort(x_{(u_{11},\ldots,u_{1r})}x_{(u_{21},\ldots,u_{2r})}).
$$
The initial polynomial of $x^u-\sort(x^u)$ is $x^u$.
\end{prop}
\begin{proof}
Since the case $r=1$ is trivial, 
assume $r\geq 2$.
Observe that by the discussion above 
for every non-face $F \subseteq rV$ of $\esd_r(\Delta^d)$
there exists a non-face $G \subseteq F$ of $\esd_r(\Delta^d)$ with $|G|=2$. 

  Since $K[\psi]\colon K[W] \to K[rV]$ takes $I'_{rT(\Delta^d)}$ to
  $I_{rT(\Delta^d)}$, it follows from Corollary \ref{inesd} that
  $\ini_{\omega'}(I'_{rT(\Delta^d)})$ is generated by square-free
  quadratic monomials of the form
  $x_{(u_{11},\ldots,u_{1r})}x_{(u_{21},\ldots,u_{2r})}$ where 
  $\{\psi(u_{11},\ldots,u_{1r}),\psi(u_{21},\ldots,u_{2r})\}$
  is {\em not} a face in $\esd_r(\Delta^d)$, that is,
  the vertices 
  $\psi(u_{11},\ldots,u_{1r})$ and $\psi(u_{21},\ldots,u_{2r})$
  are not connected by an edge in $\esd_r(T(\Delta^d))$. It follows
  that 
  $x_{(u_{11},\ldots,u_{1r})}x_{(u_{21},\ldots,u_{2r})} \in
  \ini_{\omega'}(I'_{rT(\Delta^d)})$
  if and only if there exists $i,j,k$ such that the real numbers 
  $\alpha^{ij}_k(\iota_V(\psi(u_{11},\ldots,u_{1r})))$ 
  and $\alpha^{ij}_k(\iota_V(\psi(u_{21},\ldots,u_{2r})))$ are
  non-zero and have
  opposite signs. Therefore
\begin{eqnarray*}
&&\omega'(u_{11},\ldots,u_{1r}) + \omega'(u_{21},\ldots,u_{2r})\\
&=&
\omega(\psi(u_{11},\ldots,u_{1r})) + \omega(\psi(u_{21},\ldots,u_{2r}))\\
&>&    
\omega(\psi(u_{11},\ldots,u_{1r}) + \psi(u_{21},\ldots,u_{2r})) 
\end{eqnarray*}  
On the other hand, since 
$x_{(v_{11},\ldots,v_{1r})}x_{(v_{21},\ldots,v_{2r})}
=\sort(x_{(u_{11},\ldots,u_{1r})}x_{(u_{21},\ldots,u_{2r})})$ 
is sorted we have that
$$
\omega'(v_{11},\ldots,v_{1r}) + \omega'(v_{21},\ldots,v_{2r})
=    
\omega(\psi(v_{11},\ldots,v_{1r}) + \psi(v_{21},\ldots,v_{2r})). 
$$
It follows that
$
x_{(u_{11},\ldots,u_{1r})}x_{(u_{21},\ldots,u_{2r})} 
-
x_{(v_{11},\ldots,v_{1r})}x_{(v_{21},\ldots,v_{2r})}
\in I'_{rT(\Delta^d)}$
with initial polynomial
$x_{(u_{11},\ldots,u_{1r})}x_{(u_{21},\ldots,u_{2r})}$. 
This concludes the proof.
\end{proof}
Next we deal with an arbitrary 
abstract
simplicial
complex $\Delta$.

\begin{thm}
\label{second_step}
The ideal $\ini_{\omega'}(I'_{rT(\Delta)})$ 
is generated by the initial polynomials of
the union $\Gc'$ of the binomials in $\Gc$ and
the sorted square-free monomials $x^u$ 
with $\psi(\supp(u))$ a non-face of $\esd_r(\Delta)$.
\end{thm}
\begin{proof}
The bijection $\psi \colon W \to rV$ induces an isomorphism $\psi
\colon K[W] \to K[rV]$ with $\psi(x^u) = \prod_{w \in W}
x_{\psi(w)}^{u(w)}$. 
It follows 
from \ref{toric_proj} that
\begin{eqnarray*}
&& K[W]/\ini_{\omega'} (I'_{rT(\Delta)}) \\
&\cong& 
K[rV]/
\ini_\omega (I_{rT(\Delta)})
\\
&=&
K[rV]/
\ini_\omega (I_{rT(\Delta^d)}) 
+ 
(\psi(x^u) \colon  \supp(\psi(x^u)) \not\in \Pi_{\sd_\omega rT(\Delta)})
\\
&\cong&
K[\esd_r(\Delta^d)]/
(\psi(x^u) \colon x^u \text{ square-free and sorted, } \supp(\psi(x^u))
\not\in 
\Pi_{\sd_\omega rT(\Delta)}). 
\end{eqnarray*}
Observe that for a sorted 
monomial
$x^u = \prod_{i=1}^s x_{(u_{i1},\dots,u_{ir})} \in K[W]$
we have that  $\supp(\psi(x^u)) \not\in \Pi_{\sd_\omega rT(\Delta)}$
if and only if 
$\{u_{ir}\colon 1 \leq i \leq s,\ 1\leq j \leq r \} \not\in \Delta$.
This is exactly the case if $\supp(\psi(x^u))$ is not a face of
$\esd_r(\Delta)=\Pi_{\sd_\omega rT(\Delta)}$. 
Hence
$$
K[W]/\ini_{\omega'} (I'_{rT(\Delta)})
\cong
K[\esd_r(\Delta)].
$$
\end{proof}
For the rest of this section $K$ denotes a field.
For a $\Z$-graded $K$-algebra $R$ we denote
with $H(R,n)=\dim_K R_n$ for $n \in \Z$ the
Hilbert function of $R$.
If $R={K[W]}/L$ for a graded ideal $L \subset {K[W]}$
and for a finitely generated graded module $M$
we denote with
$\projdim_R(M)=\sup\{i \in \N \colon \Tor_i^{K[W]}(R,K)\neq 0\}$
the {\em projective dimension} of $M$
and with
$\reg_R(M)=\sup\{j \in \Z \colon \Tor_i^{K[W]}(R,K)_{i+j}\neq 0\text{ for some }i \geq 0\}$
the {\em Castelnuovo-Mumford regularity} of $R$.
A standard graded $K$-algebra $R$ is said to be Koszul if
$\reg_R(K)=0$ where $K$ is regarded as an $R$-module.
For the Cohen-Macaulay and Gorenstein property of rings 
see for example \cite{BRHE}.
The next corollary lists some algebraic
consequences for the face rings.
\begin{cor}
\label{cor_esd2}
Let $\Delta$ be an abstract simplicial complex on $\{1,\dots,d+1\}$ and let $K$
be a field.
We have:
\begin{enumerate}
\item
$
\Hilb(K[\esd_r(\Delta)],n)= \Hilb(K[\Delta]^{(r)},n) \text{ for } n\geq 0.
$
\item
$K[\Delta]^{(r)}$ is 
Cohen-Macaulay, Gorenstein or Koszul if
$K[\esd_r(\Delta)]$ has one of these properties.  
\item
$\projdim_{K[W]}(K[\esd_r(\Delta)])\geq \projdim_{K[W]}(K[\Delta]^{(r)})$.
\item
$\reg_{K[W]}(K[\esd_r(\Delta)])\geq \reg_{K[W]}(K[\Delta]^{(r)})$.
\end{enumerate}
\end{cor}
\begin{proof}
Using the fact that every weight order $\omega$ can be refined to a
monomial order this follows from standard arguments.
See for example Bruns and Conca \cite{BRCO}.
\end{proof}

\begin{rem}\rm
By a theorem of
Backelin and Fr\"oberg \cite{BAFR} (see also Eisenbud, Reeves and Totaro \cite{EIRETO})
one knows that the $r$-th Veronese algebra of $K[\Delta]$
is Koszul for $r\gg 1$. One could hope that this property is
inherited for $\esd_r(\Delta)$ for $r\gg 1$. This is however not the case.
Let for example 
$\Delta$ be the set of subsets $F\subset \{1,2,3,4\}$ such that
$F\neq \{1,2,3,4\}$.
Then
$
x_{(1,4,4,\ldots)}\cdot x_{(2,4,4,\ldots)} \cdot x_{(3,4,4,\ldots)}
=0 
\text{ in } K[\esd_r(\Delta)].
$
But the elements
$
x_{(1,4,4,\ldots)}\cdot x_{(2,4,4,\ldots)},
$
$
x_{(1,4,4,\ldots)}\cdot x_{(3,4,4,\ldots)}
$
and 
$
x_{(2,4,4,\ldots)}\cdot x_{(3,4,4,\ldots)}
$
are all non-zero in $K[\esd_r(\Delta)]$.
Hence
$
x_{(1,4,4,\ldots)}\cdot x_{(2,4,4,\ldots)} \cdot x_{(3,4,4,\ldots)}
$
belongs to a minimal system of generators for the
defining ideal of 
$\esd_r(\Delta)$, but is not a quadratic monomial.
The Koszul property would imply that the defining ideal
of  $\esd_r(\Delta)$ is generated by quadratic monomials.
Hence $K[\esd_r(\Delta)]$ is not Koszul for any $r \geq 1$.
\end{rem}


\end{document}